\documentclass[preprint,review,letterpaper,11pt,numbers,sort,compress]{elsarticle}
\usepackage{amssymb}
\usepackage{amsthm}
\usepackage{amsmath}
\usepackage{bm}
\usepackage{algpseudocode}
\usepackage{algorithm}
\usepackage{multirow}
\usepackage{graphicx,subfigure}
\usepackage{caption}
\usepackage{graphics}
\usepackage{epstopdf} 
\usepackage{float}
\usepackage{eurosym}
\usepackage{anysize}
\usepackage{amsmath}
\usepackage{amsthm}
\usepackage{amsfonts}
\usepackage{fancyhdr} %%Fancy headings
\usepackage{longtable}
\usepackage{amssymb}
\usepackage{color}
\usepackage{xcolor}
\usepackage{colortbl}
\usepackage{amssymb}
\usepackage{amsmath}
\usepackage{enumerate}
\usepackage{flushend}
\usepackage{lipsum}
\usepackage{float}
\usepackage{multicol}
\usepackage{verbatim}
\usepackage{natbib}
\usepackage{hyperref}
\usepackage{tikz}% FOR AXIS
\usepackage{soul}  % TO STRIKE OFF \st{Hellow world}
\usepackage{calrsfs}
\usepackage[utf8]{inputenc}
\usepackage {csquotes}

\newcommand {\eref}{Eq.~\ref}
\newcommand {\fref}{Fig.~\ref}

\usepackage{mathtools}

\DeclareMathAlphabet{\pazocal}{OMS}{zplm}{m}{n}

\begin{document}

%%%%%%%%%%%%%%%%%%%%%%%%%%%%%%%%%%%%%%%%%%%%%%%%%%%
%%%%%%%%%%%%%%%%%%%%%%%%%%%%%%%%%%%%%%%%%%%%%%%%%%%
%%%%%%%%%%%%%%%%%%%%%%%%%%%%%%%%%%%%%%%%%%%%%%%%%%%
%%%%%%%%%%%%%%%%%%%%%%%%%%%%%%%%%%%%%%%%%%%%%%%%%%%
\begin{frontmatter}

\title{Role of coupled electrochemistry and stress on the Li-anode instability: A continuum approach}

\author{Shabnam Konica, Brian W. Sheldon, Vikas Srivastava\corref{cor1}}

\address{School of Engineering, Brown University\\
 Providence, RI, USA}
\cortext[cor1]{Corresponding author, email:Vikas$\_$Srivastava@brown.edu}

\begin{abstract}
 We present a coupled mechanistic approach that elucidates the intricate interplay between stress and electrochemistry, enabling the prediction of the onset of instabilities in Li-metal anodes and the solid electrolyte interphase (SEI) in liquid-electrolyte Li-metal batteries. Our continuum theory considers a two-way coupling between stress and electrochemistry, includes Li and electron transport through SEI, incorporates effects of Li viscoplasticity, includes SEI and electrolyte interface surface energy and evaluates crucial roles of these mechanistic effects on the continuously evolving anode surface due to the viscoplastic deformation of lithium. In the model, spatial current density evolves with the stress-induced potential across the deformed anode/SEI interface.  We assume SEI as a homogeneous, artificial layer on the Li-anode, which allows the investigation of the mechanical and electrochemical properties of the SEI systematically. Subsequently, we solve a set of coupled electrochemistry and displacement equations within the SEI and anode domains. The model is implemented numerically by writing a user element subroutine in Abaqus/Standard. We conduct numerical simulations under various galvanostatic conditions and SEI properties and predict conditions for anode instability. We find that Li viscoplasticity is one of the key attributes that drives instability in the Li-anode and show that applying a soft artificial SEI layer on the Li-anode to minimize viscoplastic deformation can be an effective method. We also report the role of artificial SEI elasticity and thickness on anode stability. Selected stability maps are provided as a design aid for artificial SEI.
\end{abstract}


\end{frontmatter}

%%%%%%%%%%%%%%%%%%%%%%%%%%%%%%%%%%%%%%%%%%%%%%%%%%%
%%%%%%%%%%%%%%%%%%%%%%%%%%%%%%%%%%%%%%%%%%%%%%%%%%%
%%%%%%%%%%%%%%%%%%%%%%%%%%%%%%%%%%%%%%%%%%%%%%%%%%%
%%%%%%%%%%%%%%%%%%%%%%%%%%%%%%%%%%%%%%%%%%%%%%%%%%%
\section{Introduction}\label{intro}
Li metal batteries have gained special attention for future electric vehicle applications due to their highest theoretical capacity  (3860mAh/g) and lowest electronegative electrochemical potential (-3.04 vs. standard hydrogen electrode) \cite{kushima2017liquid}. However, these batteries have seen limited success due to uncontrolled dendrite growth, dead lithium formation, fracture, and regrowth of SEI during cyclic charging and discharging. Although poor and inhomogeneous electrochemistry is a known factor driving anode instability, little has been known of the origin and role of the stress of the Li and SEI layers during electroplating and stripping \cite{cho2020stress}. However, the local large volume change leading to the fracture of SEI, dendritic vs. mossy Li-growth morphology, Li creep and low yield strength causing irreversible surface evolution indicate that mechanics play one of the crucial roles in the failure of Li-metal batteries \cite{eftekharnia2023understanding, wang2018stress, kushima2017liquid, goodenough2010challenges, monroe2004effect, narayan2020modeling}.  \textit{All this information requires a more in-depth computational analysis of the two-way coupling between stress and electrochemistry when lithium undergoes viscoplastic deformation and cyclic loading conditions}. This unexplored aspect is the central focus of this paper.
\\
The stability of lithium anodes depends on various factors, including the behavior of the lithium metal (viscoplasticity and creep), selection of a compatible separator (wettability, porosity, MacMullin number, etc.) \cite{eftekharnia2019toward}, and the nature of the solid electrolyte interphase \cite{liu2020review}. Achieving anode stability has been the focus of numerous experimental investigations, such as using artificial solid electrolyte interphase \cite{cho2020stress, kang2021artificial, subramanya2020carbon, choudhury2019highly}, modifying the structure or alloying of the lithium anode \cite{wan2020mechanical}, and using electrolytes and solid electrolyte interphase with better surface tension \cite{wang2017towards}. For solid electrolyte batteries, ceramics with high stiffness are often used to induce elastic stress that helps to suppress protrusion growths. However, one crucial question that remains is determining the optimal thickness of the artificial solid electrolyte interphase, which can prevent instability by providing both high surface tension and elastic stress significantly diminishing the transport properties of the solid electrolyte interphase. \par
Analytical investigations have been conducted to study anode surface instability in solid-state batteries \cite{monroe2004effect, barai2017lithium, ahmad2017stability, ahmad2017role, tikekar2016stabilizing, tikekar2016design}. Monroe and Newman's work showed that surface perturbations amplify local current density, leading to filament or dendrite growth \cite{monroe2004effect, monroe2005impact}. They suggested that a suppressor with a shear modulus twice that of Li could effectively block dendrite growth \cite{monroe2005impact}. While Stone et al. supported this concept \cite{stone2011resolution}, recent operando microscopic studies suggest otherwise, demonstrating dendrite growth even within solid ceramic electrolytes \cite{kazyak2020li, manalastas2019mechanical}. On the other hand, polymer nanocomposite separators with a modulus of approximately 1 MPa showed resilience in sustaining high cycle numbers \cite{choudhury2019highly}. One limitation in Monroe and Newman's study is the inconsistency in force equilibrium as they assumed external traction, leading to non-physical measures in their calculation \cite{deshpande2023models, mcmeeking2019metal}. \par
Ahmed and Viswanathan expanded on Monroe and Newman's stability criterion by introducing a density-driven mechanism, showing that a softer electrolyte can suppress instability if the molar volume of Li in the liquid electrolyte exceeds that in solid Li \cite{ahmad2017stability}. They later refined the criteria to consider anisotropic material properties \cite{ahmad2017role}. However, the issue of lacking force equilibrium, as noted in Monroe and Newman's work \cite{monroe2005impact}, persisted in their solution as well. \par
Barai et al. \cite{barai2017lithium, barai2016effect} developed a finite element model, building upon Monroe and Newman's analytical framework, to analyze stress and lithium transportation separately. They applied initial force equilibrium by pressing a rough Li-anode on a flat solid electrolyte and determined stress at the electrode/electrolyte interface. Utilizing these stress values, they calculated Butler-Volmer current and found that perturbations can be blocked when current density remains below a critical value \cite{barai2016effect, barai2017lithium}. Additionally, they explored the impact of external compressive stress on Li transportation properties and anode stability, hypothesizing that Li's yield stress reduces protrusion growth rate \cite{barai2018impact}. However, the \textit{decoupled stress-electrochemistry approach may be inadequate} due to the continual evolution of the anode surface during plating, particularly considering Li's viscoplastic nature.\par
Tikekar et. al. \cite{tikekar2016stabilizing} showed that the inhomogeneity in the electric field can be slowed by immobilizing a fraction of anions within the electrolyte through linear stability analysis.  In a similar study, Mcmeeking et. al \cite{mcmeeking2019metal} also performed a linear stability analysis of Li-anode with a solid electrolyte,  and through electric potential perturbation demonstrated that protrusion will always grow at a high current rate, regardless of the stiffness of the solid electrolyte. Their conclusion emphasized that perturbations with long amplitudes will invariably increase in size. However, it's worth noting that the Li-creep effect should play a role in diminishing the amplitude of roughness, particularly with a very stiff ceramic electrolyte. This aspect, as addressed in the analytical study by Zhao et al. \cite{zhao2022laplace}, deserves consideration in further analyses.
 \par
The stability analysis discussed earlier was conducted for all-solid-state batteries, where the electrolyte thickness is predetermined and thus does not warrant separate consideration. However, for artificial SEI, thickness is also one key factor. In addition, the surface tension of the liquid electrolyte/SEI interface and residual stress at the SEI layer play significant roles in the anode instability. Li-viscoplasticity coupled with SEI stress can cause wrinkling of the SEI, leading to ratcheting and delamination of the SEI upon cyclic charging and discharging \cite{liu2019wrinkling, cho2020stress}. Thus, to fully comprehend the mechanics of anode instability, all these factors demand proper attention. \par
In this work, we develop a fully coupled electro-chemo-mechanical model at a continuum level to investigate the anode and SEI instability mechanism. In this context, we are defining SEI as a thin, continuous layer of material, regardless of its origin. This includes both SEI that forms naturally through reactions, and artificially created pre-formed layers of SEI. To mimic the plating and stripping of Li-metal at the anode, we adopt a simple Li-growth model following the theory by \cite{narayan2020modeling} and consider large deformation viscoplasticity for the constitutive description of Li-metal \cite{anand2019elastic}. To predict the stability criterion of anode and SEI, we use a stress-based potential parameter as described in Monroe and Newman's work \cite{monroe2005impact}. The key contributions of the present work are summarized below:
\begin{itemize}
    \item A two-way coupling between stress and electrochemistry.
    \item Li-ion and electron transportation theory through the SEI.
    \item Incorporation of Li-viscoplasticity in the coupled electrochemistry and mechanics to investigate whether its role is stabilizing or destabilizing during cyclic plating and stripping.
    \item Incorporation of surface energy in a fully computational theory and study its effect on anode stability.
\end{itemize}
We present the following key numerical results in the context of present theory:
\begin{itemize}
    \item Mechanism of stress-coupled instability in the SEI and anode.
    \item Role of SEI modulus, thickness, surface energy and residual stress on the stability.
    \item Role of SEI inhomogeneities (voids, defects, inhomogeneous material properties
    \item Composite SEI layers and implication on Li-viscoplasticity.
\end{itemize}

The paper is structured as follows: In the first two sections, we summarize theories about Li-ion transport and stress evolution at the SEI, as well as the development of a growth model for a hypothetical interphase that resembles lithium plating and stripping. In section three, we provide the form of the equilibrium equation in the presence of surface tension. In section four, we delve into the strategic implementation of the Butler-Volmer current on a finite surface. The governing partial differential equations, along with corresponding initial and boundary conditions, are explained in section five. In section six, we present the numerical simulation results, featuring stability maps of the anode surface as a function of SEI properties. Finally, in section seven, we conclude with closing remarks.

%%%%%%%%%%%%%%%%%%%%%%%%%%%%%%%%%%%%%%%%%%%%%%%
\section{Theory of Li-anode and SEI: Ion transportation and constitutive behavior }
\begin{figure}[h!]
  \centering
  \begin{tabular}{c}
     \includegraphics[width=0.9\textwidth]{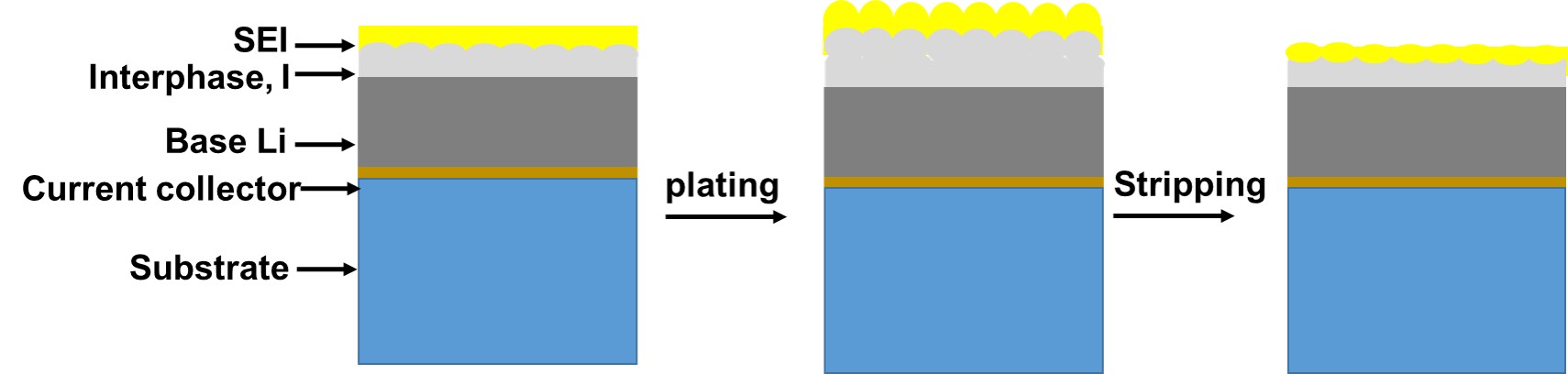}
     \centering
     \end{tabular}   
  \caption{Finite element model setup for Li-anode and SEI: Li is plated and stripped at the interphase I}
    \label{FEA_model_1}
\end{figure}

We limit our investigations to predict the onset of instability at the anode/SEI due to coupled electro-chemo-mechanics, where we solely interest ourselves in the amplification of surface perturbations, and assume that the occurrence of fracture of the SEI and dendrite growth happens after instability. \fref{FEA_model_1} shows the block diagram of the anode and SEI instability problem in a symmetric Li half cell. The half cell is constructed of a thin layer of SEI, followed by an ultra-thin interphase layer I, a thick layer of an initial Li layer (base Li), and a current collector on a thick quartz substrate. The Li-ion migration within the SEI layer is solved in a coupled electrochemistry and large deformation problem and a reaction current density is computed assuming Butler-Volmer kinetics at the SEI/I interface, which depends on the voltage, Li-ion concentration, and stress at the interface. We then solve the Li-migration through the interphase I using a Fickian-type diffusion with a very high diffusion coefficient where the reaction current density acts as a source of Li-ions. The flow of Li-ions causes the interphase to swell, which phenomenologically represents the plating of Li at the anode. We replaced the sharp Li/SEI interphase with an ultra-thin metallic layer denoted as $\textrm I$ and used bulk growth of the interphase I to model the surface growth process, similar to \cite{narayan2020modeling}. The interphase $\textrm I$ possesses the same material properties as Li metal and can grow or shrink dynamically by absorbing or exuding Li-ions during charging or discharging, respectively, to emulate the plating and stripping process. Finally, Li-ions migration causes the viscoplastic deformation of the interphase. 
%%%%%%%%%%%%%%%%%%%%%%%%%%%%%%%%%%%%%%%%%%%%%%
%%%%%%%%%%%%%%%%%%

\subsection{Theory of Ion transportation through SEI}
In this section, we provide a summary of Li-ion and electron transportation through SEI \cite{narayan2021coupled, ganser2019extended, ganser2019finite, bistri2023continuum, zhao2022phase}. 
\subsubsection{Balance of Li-ion and electron}
Assuming that only Li ions flow through the SEI layer and anions are immobile, the local Li-ion balance equation is expressed as, 
\begin{equation}
    \dot c_{R}=-\textrm{Div} \mathbf j_R
    \label{cRbalance}
\end{equation}
where $c_R$ is the local ion concentration at the SEI layer (mol of cation per unit volume) and $\mathbf j_R$ is the referential flux (per unit area per unit time) that depends on the gradient of electrochemical potential and Li-ion mobility. The electrochemical potential is defined as,
\begin{equation}
    E_{EC}=\mu_0+ \textrm {ln} \left( \frac{\bar c}{(1-\bar c)} \right)+ F\phi \label{specific_pot} 
\end{equation}
where $\bar c$ is the normalized Li-ion concentration defined as, $\bar c=c_R/c_{\textrm{max}}$, $\textrm F$ is the Faraday constant, $\phi$ is the voltage. Next, mobility can be defined as, $ \mathbf M=\frac{D}{R\vartheta}c_{R}\left (1-\frac{c_{R}}{c_{max}}\right)\mathbf{I} $, where,  $D$ is the diffusion coefficient of Li ions inside SEI, $R$ is the universal gas constant, and $\vartheta$ is the operating temperature. Using the fact that, $ \frac{\partial\mu_{chem}}{\partial c_R}= \frac{RT}{c_{R}(1-\bar c)}$; $\frac{\partial\mu_{mech}}{\partial p_{\textrm{SEI}}} =-\Omega_{SEI}$ where, $p_{\textrm{SEI}}$ is the hydrostatic stress and $\Omega_{\textrm{SEI}}$ is the molar volume of Li-ion through SEI,  using \eref{cRbalance}, we can rewrite the  Nernst-Planck equation for Li-ion migration through SEI under a given electric field as,
\begin{equation}
    \dot c_R=-\textrm{Div}D\nabla c_R+\frac{Dc_{R}}{R\vartheta}\left (1-\frac{c_{R}}{c_{max}}\right)\left( F\nabla \phi + p \Omega_{SEI} \right)
\end{equation}
The balance equation of electrons for SEI considering charge neutrality is described through Poisson's equation of voltage \cite{chen2015modulation, bistri2023continuum, narayan2021coupled, barai2017lithium} as,
\begin{equation}
    \textrm{Div} (\kappa \nabla \phi)=0 \label{voltagePDE}
\end{equation}
where $\kappa$ is the electric conductivity of SEI.
\subsubsection{Electrochemical interfacial kinetics}

During the plating process, Li-ion undergoes an electrochemical ion-transfer reaction, where it takes an electron and deposits as Li-metal at the Li-anode. This phenomenon is quantified by the Butler-Volmer equation \cite{deshpande2023models}, a pivotal tool expressing the rate at which ions traverse the interface per unit electrode surface area. The equation introduces an exchange current density $i_0$, a Bronsted–Evans–Polyani symmetry factor $\beta$ $(0\leq \beta\leq 1)$ \cite{ganser2019extended, newman2004electrochemical}, and the surface overpotential $\eta_s$, defined as: 
\begin{equation}
    i=i_0 \; \left[\textrm{exp}\left [ \frac{\beta F \eta_s}{RT}\right]-\textrm{exp}\left[- \frac{(1-\beta)F \eta_s}{RT}\right] \right] \label{BV1}
\end{equation}
Here, $i_0$ is the exchange current density of Li ions moving from electrolyte to anode, $\beta$ is Bronsted–Evans–Polyani symmetrey factor $(0\leq \beta\leq 1)$ \cite{ganser2019extended, newman2004electrochemical} and $\eta_s$ is the surface overpotential and defined as,
\begin{equation}
    \eta_s=\phi_{\textrm{electrode}}-\phi_{\textrm{electrolyte}}-\phi_{\textrm{eq}}
\end{equation}
Here $\phi_{\textrm{electrode}}$ is the electric potential at the Li-anode and $\phi_{\textrm{electrolyte}}$ is the electrolyte electric potential. The equilibrium potential can be computed from the energy barrier arising from the chemical potential difference between the anode and electrolyte and the mechanical stress difference at the interface. \par
Numerous studies, such as \cite{ganser2019extended, mcmeeking2019metal, monroe2004effect, bazant2013theory}, highlight the dependence of the exchange current density on the equilibrium potential. Bazant and colleagues \cite{bazant2013theory}, adopting Bronsted–Evans–Polyani theory, showcased that configurational entropy's impact on the energy barrier is negligible. Monroe and Newman extended this by introducing an additional symmetry factor $\alpha$ to consider stress effects, leading to the expression for exchange current density defined as,
\begin{equation}
    i_0=\Tilde{i_0} \textrm{exp}\left[\frac{[\alpha-(1-\beta)]\;(\mu^E_{\textrm{electrode}}-\mu^E_{\textrm{electrolyte}})}{RT} \right] \label{i0Mon}
\end{equation}
Here, $\mu^E_{\textrm{electrode}}$ and $\mu^E_{\textrm{electrolyte}}$ denote the electric potential at the deformed electrode and electrolyte interface, respectively. Notably, when $\alpha=1-\beta$, the stress contribution from \eref{i0Mon} vanishes. On the other hand, setting $\alpha=1$ gives a stress effected exchange current density expression used in \cite{monroe2004effect, monroe2005impact} and expressed by,
\begin{equation}
     i_0=\Tilde{i_0} \textrm{exp}\left[\frac{\beta\;(\mu^E_{\textrm{electrode}}-\mu^E_{\textrm{electrolyte}})}{RT} \right] \label{iO-2}
\end{equation}
Monroe and Newman \cite{monroe2004effect} demonstrated that in the absence of Li-ion influence, the quantity $\mu^E_{\textrm{electrode}}-\mu^E_{\textrm{electrolyte}}$ is equivalent to the mechanical pressure-induced potential difference between the electrode and the electrolyte layer, denoted as $(p_{\textrm{electrode}}\Omega_{\textrm{electrode}}-p_{\textrm{electrolyte}}\Omega_{\textrm{electrolyte}})$, where the subscripts indicate the the location and $p$ and $\Omega$ denote pressure and molar volume, respectively.  A specific form of this quantity in terms of external traction is also present in their follow-up work \cite{monroe2005impact}. We adopt a similar expression to compute the exchange current density at the Li-anode and SEI interface. Utilizing \eref{iO-2} in \eref{BV1}, we arrive at the complete Butler-Volmer reaction current expression at the anode/SEI interface:
\begin{equation} \label{BVcurrent}
     i=\Tilde{i_0} \textrm{exp}\left[\frac{\beta\;(\mu^E_{\textrm{electrode}}-\mu^E_{\textrm{electrolyte}})}{RT} \right]\; \left[\textrm{exp}\left [ \frac{\beta F \eta_s}{RT}\right]-\textrm{exp}\left[- \frac{(1-\beta)F \eta_s}{RT}\right] \right]
\end{equation}
It's noteworthy to mention that Ganser et al. introduced a new factor, $\delta$, to account for the effect of equilibrium potential on exchange current density. However, we neglect the influence of chemical concentration differences while computing the reaction current due to the negligible thickness of the SEI compared to the anode.
\subsubsection{Specific form of free energy, and Cauchy stress}
The total free energy contribution comes from mechanical energy and energy due to the diffusion of Li-ions into SEI. Thus, we write the following separable form of the free energy:
\begin{equation} \label{freeEnergy}
    \psi_R=\mathit J^s \;\left( \frac{1}{2}\mathbf E^e:\mathbb C\mathbf E^e \right) +R\vartheta c_\textrm{Rmax}[\bar c\textrm{ln}\bar c+(1-\bar c)\textrm{ln}(1-\bar c)],   
\end{equation}
where $\mathit J^s=\mathit J^s= 1+\Omega_{\textrm{SEI}}(c_R-c_{R0})$ is the volume change due to swelling of SEI, $\mathbb C=2G\mathbb I+(K-\frac{2}{3}G)\mathbf I\otimes \mathbf I$ is the elasticity tensor and $\mathbb I$ and $\mathbf I$ are fourth-order and second-order Identity tensor, respectively. $\mathbf E^e$ is the elastic logarithmic tensor. Starin takes an additive decomposition into the elastic and swelling components as, 
\begin{equation}
    \mathbf E=\mathbf E^e+\frac{1}{3}\textrm{ln}(1+\Omega_{\textrm{SEI}}(c_R-c_{R0})) \textbf I
\end{equation}
%
%\paragraph{\textbf{cauchy Stress}}
The Cauchy stress can be defined by the equation as follows:
\begin{equation}\label{cauchyfinal}
    \mathbf T=\mathit J^{e-1}\left[ 2G\mathbf E_H^e+\left( K-\frac{2}{3}G\right)(\textrm{tr}\mathbf E^e_H) \mathbf I\right]-K\textrm{ln}\big[1+\Omega_{\textrm{SEI}}(c_R-c_{R0})\mathbf I\big]
\end{equation}
where $\mathbf E_H^e$ is the spatial elastic logarithmic strain \cite{narayan2021coupled}, $G$ is the shear modulus and $K$ is the bulk modulus of SEI. $\mathit J^e$ is the determinant of the elastic deformation gradient and $\mathbf I$ is the identity tensor.
%%%%%%%%%%%%%%%%%%%%%%%%%%%

\subsection{Constitutive theory and growth of the interphase \textit I} \label{interphaseCons}
\subsubsection{Li-ion balance at the interphase I due to Li-influx}
We hypothesize that the influx of Li-ions into interphase \textit I is contingent upon the reaction current density $i$ at the SEI/I interface. In our assumption, the rapid diffusion and high conductivity within interphase \textit I $(\mathbf M_I >> \mathbf M_{SEI})$ ensure that the ion concentration gradient remains practically close to zero. This not only aligns the rate of Li-ion concentration within the interphase with the reaction current density but also alleviates the numerical challenge of assigning the reaction current density at each Gauss point within \textit I. Our formulation of the local mass balance at interphase I is expressed as:
\begin{equation}
    \dot c_{R}=\textrm{Div}\mathbf j_R, \textrm{  with,   } \; \; \; \; \; \mathbf j_R=D_I\nabla c_R+\frac{D c_R}{R\vartheta}p_I\Omega_I  
\end{equation}
Here, $D_I$ is a superficial high diffusivity value that we choose in a way that there is no concentration gradient between the SEI/I interface to the I/Base Li interface, $p_I$ is the pressure at interphase \textit I, and $\Omega_I$ is the molar volume of interphase \textit{I}. 

%%%%%%%%%%%%%%%%%

\subsubsection{Kinematics of interphase I}
Our interphase theory is built upon Narayan and Anand's model \cite{anand2019elastic, narayan2020modeling} regarding Li-viscoplasticity and growth. Herein, we provide a concise overview of the constitutive theory. The theory encompasses a multiplicative decomposition of the deformation gradient, an anisotropic growth evolution for the interphase, and viscoplastic deformation \cite{kothari2019thermo, rejovitzky2015theory, di2015diffusion, niu2023modeling, Srivastava2011STRESSPIPE, bai2021continuum}. The multiplicative decomposition of the deformation gradient is expressed as:
\begin{equation}
    \mathbf F=\mathbf F^e \mathbf F^p \mathbf F^g, \textrm{ with } \textit J^g=1+\Omega_I c_R; \textrm { and   } \dot {\mathbf F}^g=\mathbf D^g\mathbf F^g; \text{  \; \; } \dot {\mathbf F}^p=\mathbf D^p\mathbf F^p
\end{equation}
with $\mathit J^g= \textrm{det} \mathbf F^g$, and $\mathbf F^g$ and $\mathbf F^p$ denote the distortion due to growth \cite{narayan2020modeling,Zhong2021AGrowth} and plasticity, respectively, and $\mathbf F^e$ is the elastic deformation gradient, $\mathbf D^g$ and $\mathbf D^p$ denote the strain tensors due to growth and plastic deformation, respectively. Finally, $\Omega_I$ represents the molar volume of \textrm I.

The anisotropic plating growth of the interphase can be expressed as,
\begin{align*}
    \mathbf D^g=\dot\epsilon^g \mathbf S^g \textrm{,  with } \mathbf S^g=\sum_{i=1}^3\alpha_i \mathbf m_i\otimes \mathbf m_i; \\
    \textrm { \; where \; }  \dot\epsilon^g=\frac{\Omega_I\dot c_R}{1+\Omega \dot c_R}  \textrm{,\; 
 with \; }
    \mathbf m_i=\frac{(\mathbf F^e \mathbf F^p)^{-1}\mathbf n(\mathbf x)}{|(\mathbf F^e \mathbf F^p)^{-1}\mathbf n(\mathbf x)|}
\end{align*}
Here $\epsilon^g$ is the volumetric growth strain, $\Omega$ represents the molar volume of \textrm{I}, and $m_i$ is the growth direction. Additionally, $\mathbf S^g$ is an anisotropic growth tensor containing intermediate growth directional unit vectors,k with $\alpha_i$ as scalar constants representing growth weights and $\mathbf n(\mathbf x)$ representing the normal to the interface. \par
To solve the evolution of $\mathbf F^g$, one needs to provide an initial condition $\mathbf F^g(\mathbf X,0)=\mathbf I$.
Next, the evolution of plastic deformation gradient $\mathbf F^p$ is given as,
\begin{equation}
    \dot {\mathbf F}^p=\mathbf D^p\mathbf F^p; \textrm {  where  }  \mathbf D^p=\frac{3}{2}\bar{\epsilon}^p \big(\frac{\mathbf M^e_0}{\bar \sigma} \big); \textrm { and } \dot{\bar{\epsilon}}^p=A\textrm{exp}\left(-\frac{Q}{RT}\right)\;\left( \frac{\bar{\sigma}}{S}\right)^{1/m};
\end{equation}

Here, $\bar{\epsilon}^p$ is the equivalent plastic strain, and $\mathbf M^e_0$ is the deviatoric part of the Mandel stress, $A$ is a pre-exponential factor, $Q$ is the activation energy, $\bar{\sigma}$ is the equivalent tensile stress and $m$ is a strain-rate sensitivity parameter, $S$ is the strain-rate dependent flow resistance given as,
\begin{equation}
    \dot S=\bigg[ H_0\bigg|1-\frac{S}{S_s}\bigg|^a \;\textrm{sign}\bigg( 1-\frac{S}{S_s}\bigg)\bigg] \textrm{, with  } S_s=S_*\left[ \frac{\dot{\bar{\epsilon}}^p}{A\textrm{exp}\bigg(-\frac{Q}{RT}\bigg)}\right]
\end{equation}
Here, $S_s$ is the saturation value of flow strength at a given strain or temperature rate, and \{$H_0, S_*, a, n$\} are strain-hardening parameters.

\subsubsection{Free  energy and Stress}
The mechanical free energy takes the form,
\begin{equation}
    \psi_M=\mathit J^g \left(G|\mathbf E_0^e|^2+\frac{1}{2} K |\textrm{tr}\mathbf E_0^e|^2 \right)
\end{equation}
here $G$ and $K$ are shear and bulk modulus of Li, $\mathbf E^e_0$ is the deviatoric part of logarithmic elastic strain $\mathbf E^e$, and $\mathit J^g=\textrm{det} \mathbf F^g$. 

The equivalent tensile stress that governs viscoplastic flow is defined from the deviatoric part of Mandel stress \cite{konale2023large,Srivastava2010SMP} as,
\begin{equation}
    \bar\sigma=\sqrt{3/2}|\mathbf M^e_0|
\end{equation}
Finally, The Cauchy stress, $\mathbf T$ is given by,
\begin{equation}
  \mathbf T =\mathit J^{e-1}\left[ 2G\mathbf E_H^e+\left( K-\frac{2}{3}G\right)(\textrm{tr}\mathbf E^e_H) \mathbf I\right]
\end{equation}

Finally, for the base Li metal, we consider a large deformation elastic-viscoplastic model similar to the theory mentioned for interphase I, but without growth, such that, $\mathit J^g=1$.
\section{Role of surface tension on the equilibrium equation}

Let us consider an assembly of SEI and I in \fref{FEA_model_1}. Surface tension works at the SEI/I interface and at the SEI and liquid electrolyte interface. \par
For the SEI/I interface, the force balance requires,
\begin{equation}\label{traction1}
    \mathbf T_{\textrm{I}}\mathbf n_1=\mathbf t_1+\mathbf t_1
\end{equation}
Here $\mathbf t_1$ is the traction imparted by surface tension $\gamma_{SS}$, while $\mathbf t_2$ is the traction that comes from the load applied by SEI onto the interphase I.
Thus $\mathbf t_1=\mathbf T_{\textrm{SEI}}\;\mathbf n_1$, and $\mathbf t_2=-2\mathcal H_1 \;\gamma_{SS}$, where $\mathcal H_1$ is the mean curvature. Substituting these values into \eref{traction1}, we obtain the equilibrium force working at the SEI/I interface as,
\begin{equation}\label{force_bal1} 
    \mathbf T_{\textrm{I}}\mathbf n_1 -\mathbf T_{\textrm{SEI}}\mathbf n_1+2\mathcal H_1 \;\gamma_{SS}=0
\end{equation}
Similarly, at the SEI liquid interphase, the force balance is,
\begin{equation}
    \mathbf T_{\textrm{SEI}}\mathbf n_2=-2\mathcal H_2\;\gamma_{SL}
\end{equation},
where $\mathcal H_2$ is the mean curvature at the SEI/liquid interface.
In our simulations, we apply the normal tractions coming from the surface tensions at the respective surfaces to the load integration points on the free surface of the finite element model \cite{henann2010surface, henann2014modeling, wang2016modeling}. 
%%%%
\section{Numerical calculation of Butler-Volmer current}
Since SEI and I are discontinuous continuum bodies, one particularly challenging aspect is to compute the stress-induced chemical potential difference ($p_I\Omega_I-p_{\textrm{SEI}}\Omega_{\textrm{SEI}}$)at the SEI/I interface in a coupled finite element computational platform. Barai et al. \cite{barai2017lithium} tackled this complexity by employing a decoupled approach, applying external pressure to the anode and separator assembly to compute the force equilibrium and subsequently determine the exchange current density based on the stress calculated in the previous step. While this method works for initial steps, the \textit{cyclic nature of battery charging and discharging introduces viscoplastic deformation in the anode, causing continuous evolution of the anode/SEI surface}. Thus decoupling stress and electrochemistry will not provide an accurate representation of stress and reaction current densities\par
To address this issue, we have considered a thin surface situated between SEI and I. This surface comprises two layers of elements. The top layer of the surface shares the same properties as SEI, while the bottom layer has the same properties as interphase I. Both layers have the same nodes at the middle line or the SEI/I interface, thus making it possible to resolve the stress-induced potential at these nodes by creating an additional differential equation written as,
\begin{equation}\label{musigma}
    \mu_{\sigma}=p_I\Omega_I-p_{\textrm{SEI}}\Omega_{\textrm{SEI}}
\end{equation}
,where $\mu_{\sigma}$ is the stress-potential difference at the SEI/I interface.
This differential equation can be solved at the nodes on the SEI/I interface, enabling stress-induced potential tracking at the interface.
The surface overpotential is computed as,
\begin{equation}\label{overp}
    \eta_S=\phi_I-\phi_{\textrm{SEI}}-\Delta\phi_{eq}-\mu_{\sigma}
\end{equation}
We assume that the amount of Li-ion entering the interphase I is equal to the reaction current density at the SEI/I interface. Therefore, there is no concentration difference at the SEI/I interface. Consequently, the equilibrium potential is dependent only on the stress at the SEI/I interface. Finally, \eref{musigma} and \ref{overp} are inserted into \eref{BVcurrent} to compute the Butler-Volmer type reaction current density at the interface.

\section{Governing differential equations and the boundary conditions}

We solve a fully coupled electro-chemo-mechanical problem for SEI by solving 3 partial differential equations:
\begin{equation}
\begin{cases}
\textrm{Div}\mathbf T+\mathbf{b}=\mathbf{0} \\
\dot{c}_R=-\text{Div}\mathbf j_R \\
\textrm{Div}(\kappa\nabla\phi)=0
\end{cases}
\label{pdes}
\end{equation}
It is to be noted that, we solve all three equations mentioned in \eref{pdes} in the case of SEI, where for for the hypothetical interphase I, we solve only the momentum and Li-ion balance equation. For base Li, only the momentum equation needs to be solved. The system of equations can be solved by considering appropriate initial and boundary conditions.  Here the system of equations is solved numerically for each element by writing a user element subroutine in Abaqus/Standard (2020) \cite{AbqStan}.

\section{Numerical Simulations, Results and Discussions}
Here we present numerical examples to explore the main factors that cause instability at the SEI and the anode. We investigate the impact of three key factors: i) the amplitude and wavelength of the initial perturbation at the anode; ii) the presence of voids at the SEI/I interface; and iii) SEI single or hybrid layer nature, and its electrochemical and mechanical properties. Based on a comprehensive study of these three factors, we suggest an optimized design of SEI that eliminates the mechanical factors responsible for instability propagation, up to some safe current density. The material parameters for Li and hypothetical interphase \textit{I} are given in Table~\ref{table:MatParams}, and the diffusivity and conductivity of the hypothetical interphase \textit{I} are listed in Table~\ref{table:Kinpar}. We vary the mechanical and electrochemical properties of the SEI layer to better understand their role in the instability and report our findings in each numerical simulation section.
{\renewcommand{\arraystretch}{1.1}
\begin{table}[!h]
\caption{Material parameters of Li and hypothetical interphase \textit{I} taken from (\cite{anand2019elastic})}
\begin{center}
\begin{tabular}{|cccccccccc|}
\hline
\hline
E(GPa)  & $\nu $ & A($s^{-1})$ &  Q($\textrm{kJ}\textrm{mol}^{-1}$) & m & $S_0$ (MPa) & $S_*$(MPa) & $H_0$(MPA) & a & n \\
\hline
7.81 & 0.38 & $4.25\times 10^{4}$ & 37 & 0.15 & 0.95 & 10.0 & 2.0 & 2.0 & 0.05\\
\hline
\end{tabular}
\end{center}
\label{table:MatParams}
\end{table}}
%%%%%%%%%%%%%%%%%%%%
{\renewcommand{\arraystretch}{1.1}
\begin{table}[!h]
\caption{Numerical kinetic parameters for Interphase I}
\begin{center}
\begin{tabular}{|ccc|}
\hline
\hline
Diffusivity($m/s^2$)  & Conductivity($S/cm$) & $c_{\textrm{max}}, \textrm{mol}/{m}^3$ \\
\hline
$1\times 10^{-9}$ & 10 & $10\times 10^4$\\
\hline
\end{tabular}
\end{center}
\label{table:Kinpar}
\end{table}}
%%%%
\subsection{Modeling instability propagation in a SEI layer: Role of  initial imperfection}

\begin{figure}[h!]
  \centering
  \begin{tabular}{c}
     \includegraphics[width=0.85\textwidth]{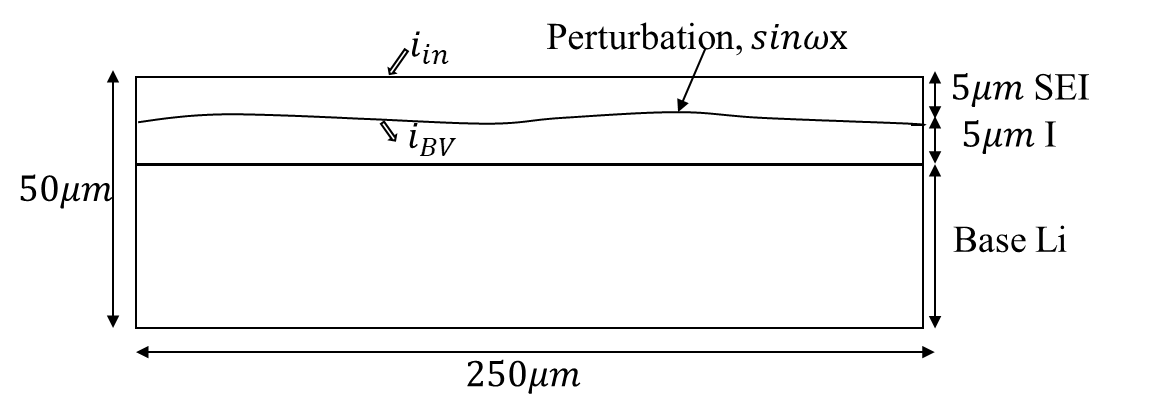}\\
     \centering
     \end{tabular}   
  \caption{Geometry of a Li anode with a 5$\mu m$ thick SEI}
    \label{5umSEI_geo}
\end{figure}

Let us consider the geometry of a Li-anode with SEI assembly as shown in \fref{5umSEI_geo}. The dimension of the specimen is 250X50 $\mu m$. A $5~\mu m$ thick SEI is considered on top of a $5~\mu m$ hypothetical interphase layer \textit{I}. There is a $0.2~\mu m$ thin copper current collector layer at the bottom of the specimen and rest $43.8 ~\mu m$ is the base Li. We consider a symmetric boundary on the left and right side of the specimen, and bottom is fixed. An initial sinusoidal perturbation of $\textrm{a\;sin}\omega \textrm {x}$ is applied at the interface of SEI and \textit I, x is the horizontal axial coordinate value. \par
%%%%% result figure 1
%%%% 
For this numerical study, we start the investigation at a perturbation wavelength of $\omega=10/\mu m$ and gradually decrease the wavelength to $500/\mu m$ and observe the anode surface condition (flat vs. wavy). \fref{5umSEI_geo} shows the finite element geometry of two such cases, when $\omega=10/\mu m$ and $100/\mu m$, respectively. The amplitude of perturbation is $a=0.4\; \mu m$ for both cases. The properties of SEI are as follows: $\textrm E=340~\textrm{MPa}; ~\nu=0.33;~ \kappa_{\textrm{cond}}=1\times 10^{-6}~S/cm;~ \textrm D=1\times 10^{-8} \textrm{cm}^2/s$, and $\Omega_{\textrm{SEI}}=2\times 10^{4} m^3/\textrm{mol}$. In the next paragraph, we will discuss the relation between initial perturbation vs. plating current on anode instability. \par
\begin{figure}[h!]
  \centering
  \begin{tabular}{c}
     \includegraphics[width=0.95\textwidth]{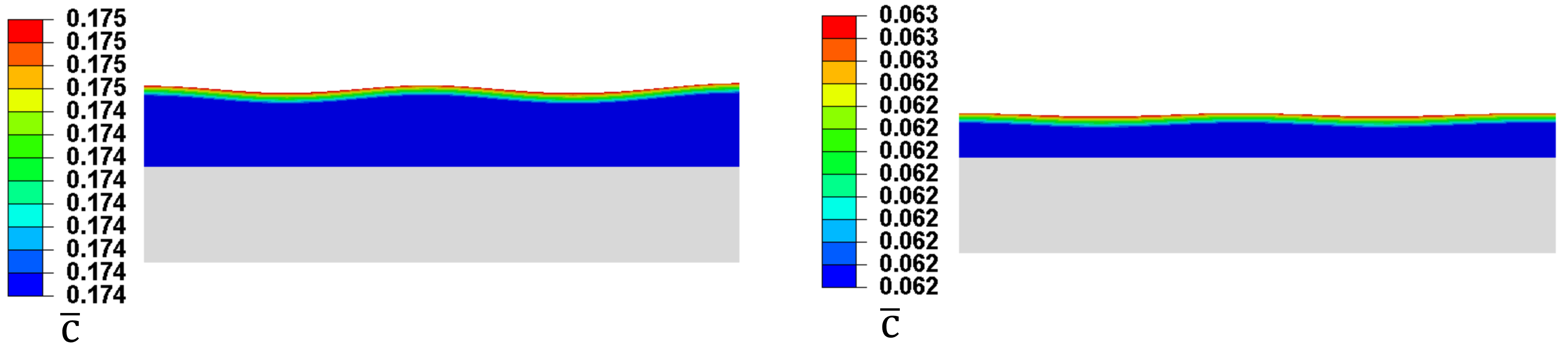}\\
     (a) \\
     \includegraphics[width=0.95\textwidth]{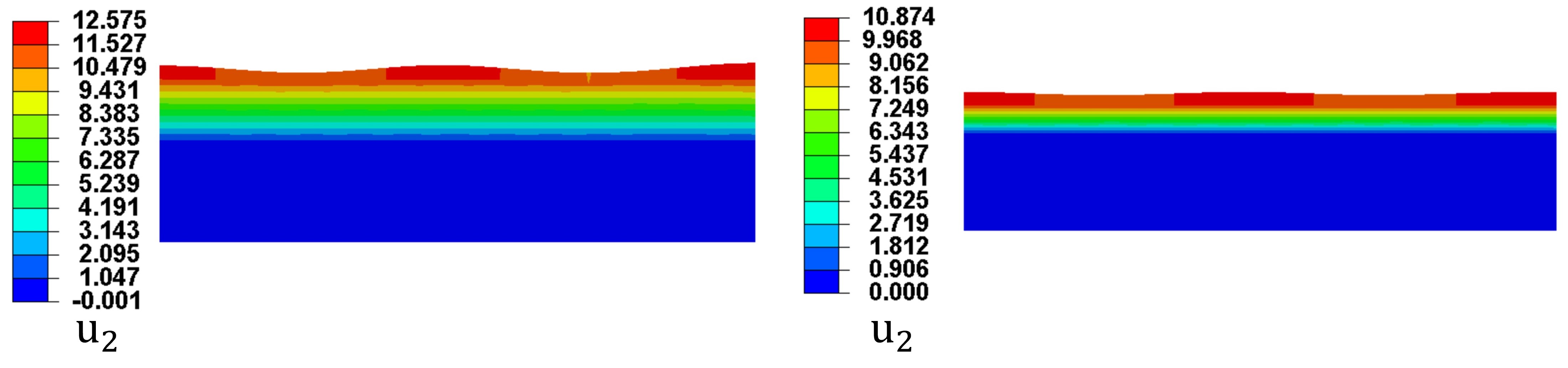}\\
     (b)\\
     \includegraphics[width=0.95\textwidth]{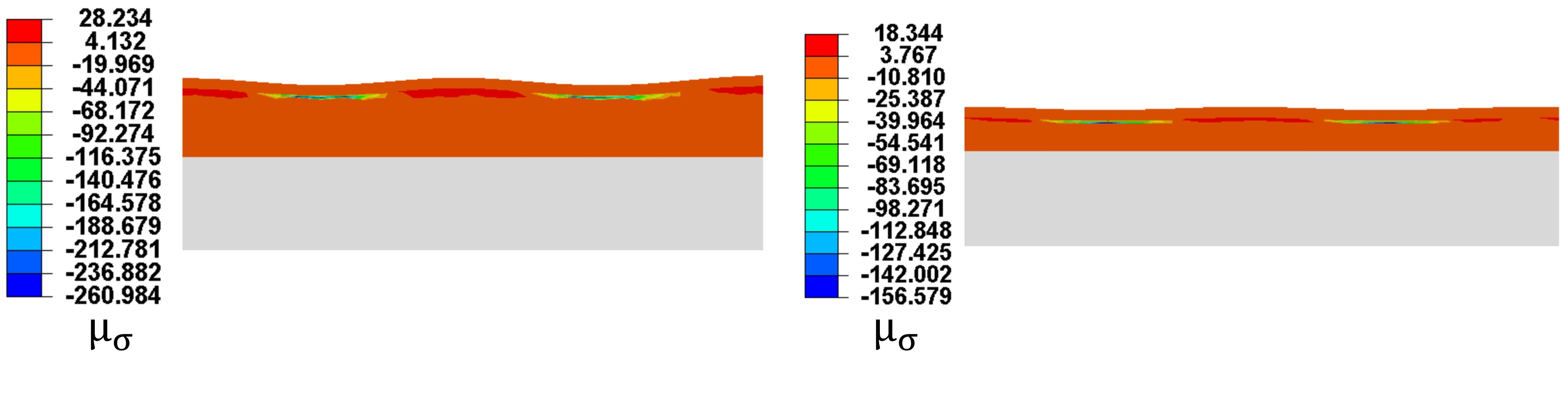}\\
     (c)
     \end{tabular}   
  \caption{Contour plots of (a) Normalized Li-ion concentration; (b) displacement in $\mu m$; (c) Stress-induced potential in J/mol plated at (left) $2.5~\textrm{mA}/{\textrm{cm}}^2$ for 1 hr; (right) plated at $0.5~\textrm{mA}/{\textrm{cm}}^2$ for 1 hr }
    \label{5um_2vs0p5_v1}
\end{figure}
\begin{figure}[h!]
  \centering
  \begin{tabular}{c}
      \includegraphics[width=0.95\textwidth]{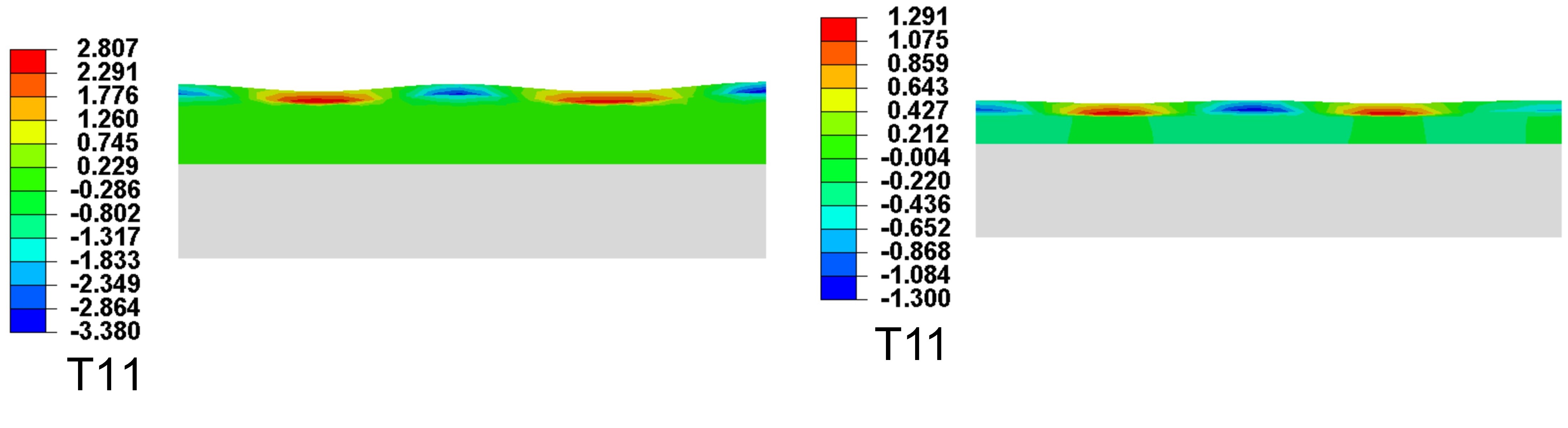}\\
     (a)\\
     \includegraphics[width=0.95\textwidth]{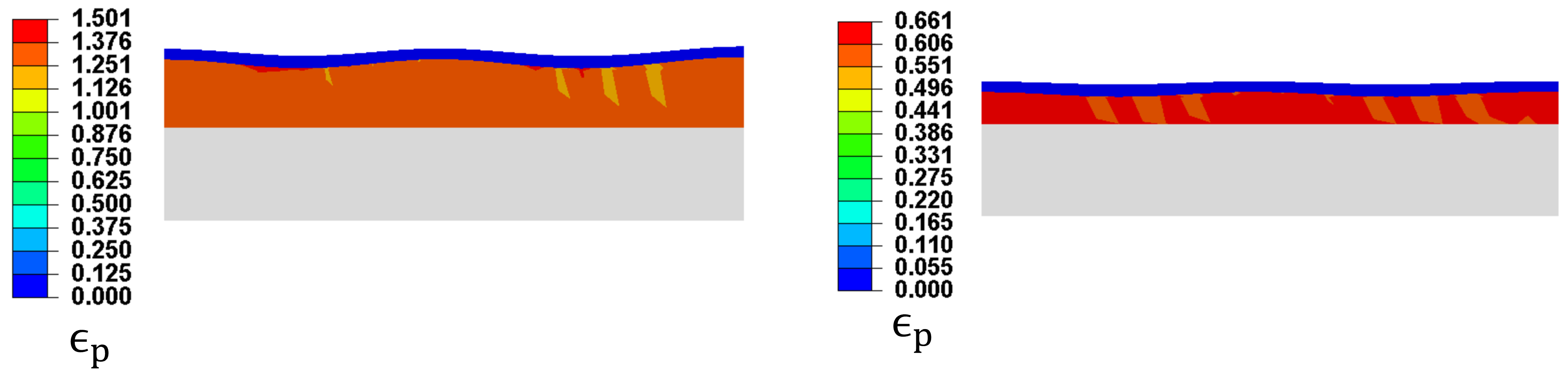}\\
     (b)
     \end{tabular}   
  \caption{Contour plots of (a) Stress T11 in MPa; and (b) Plastic strain (left) plated at $2.5~\textrm{mA}/{\textrm{cm}}^2$ for 1 hr; (right) plated at $0.5~\textrm{mA}/{\textrm{cm}}^2$ for 1 hr }
    \label{5um_2vs0p5_v2}
\end{figure}

\paragraph{\textbf{Instability propagation vs. current density}}
In the case of perturbation at $\omega=10/\mu m$, as shown in Figure \ref{5umSEI_geo}(a), we conducted several tests where we plated Li for one hour at different current densities of $[2.5, 2.0, 1.5, 1, 0.5]\textrm{ mA}/{\textrm{cm}}^2$. Our findings indicate that the anode and SEI became unstable in all test scenarios, even at a very low current density of $0.5~\textrm{mA}/{\textrm{cm}}^2$. This was demonstrated by the appearance of wrinkles on the SEI top surface, which indicates that the perturbation wavelength is too large in this case. As a representative case, we here show the results at two current densities, $2.5~\textrm{ mA}/{\textrm{cm}}^2$ and $0.5~\textrm{ mA}/{\textrm{cm}}^2$, as shown in \fref{5um_2vs0p5_v1} and \ref{5um_2vs0p5_v2}, respectively, where we show the contour plots of normalized current density, vertical displacement at anode, stress-induced potential, stress component T11 and plastic strain after 1 hour of plating.

%%%% result figure 2
 Both current rates show homogeneous Li-ion concentrations at the SEI/I interface, as seen in \fref{5um_2vs0p5_v1}(a). The higher current rate ($2.5~\textrm{mA}/{\textrm{cm}}^2$) results in a higher ion concentration, as shown by the left side contour plot of \fref{5um_2vs0p5_v1}(a).
Comparing the vertical displacement contours for the two cases from \fref{5um_2vs0p5_v1}(b), it is observed that the peaks grew significantly more than the plateau in the case of the current rate $2.5~\textrm{mA}/{\textrm{cm}}^2$. This occurs due to the higher difference in stress-induced potential $\mu_{\sigma}$ between peak and valley in case of $2.5~\textrm{mA}/{\textrm{cm}}^2$ current rate, as demonstrated in \fref{5um_2vs0p5_v1}(c) (left). Negative values at the valley reduced the exchange current density, causing the valleys to grow at a slower rate. In case of the lower current density of ($0.5~\textrm{mA}/{\textrm{cm}}^2$), almost uniform $\mu_{\sigma}$ was observed in both peaks and valleys, resulting in both peak and valley growing at almost equal rates \fref{5um_2vs0p5_v1}(b-c)(right).
\par
The axial stress contours in both cases showed a negative bending stress at the peak followed by a positive bending stress at the valley \fref{5um_2vs0p5_v2}(a). Although the same geometry and material properties were considered for both cases, the bending stress magnitudes were much higher at the higher current rate due to higher plastic strain. The plastic strain contours in \fref{5um_2vs0p5_v2}(v) showed that higher plastic strain caused higher equivalent tensile stress values at the interphase I, thus resulting in higher stress experienced by the SEI in case of $2.5~\textrm{mA}/{\textrm{cm}}^2$. \par   
\paragraph{\textbf{Instability propagation vs. perturbation wavelength}}
For a second study, we show the results with three perturbation wavelengths: $\omega=10/ \mu m,\; 100/ \mu m\; \textrm{ and } 500/ \mu m$. In this context, we examine a plating scenario at a current rate of $0.5~\textrm{mA}/{\textrm{cm}}^2$ for 1 hour, maintaining the same SEI geometry and properties as in the previous study with $\omega=10/\mu m$. \par
\begin{figure}[h!]
  \centering
  \begin{tabular}{cc}
     \includegraphics[width=0.95\textwidth]{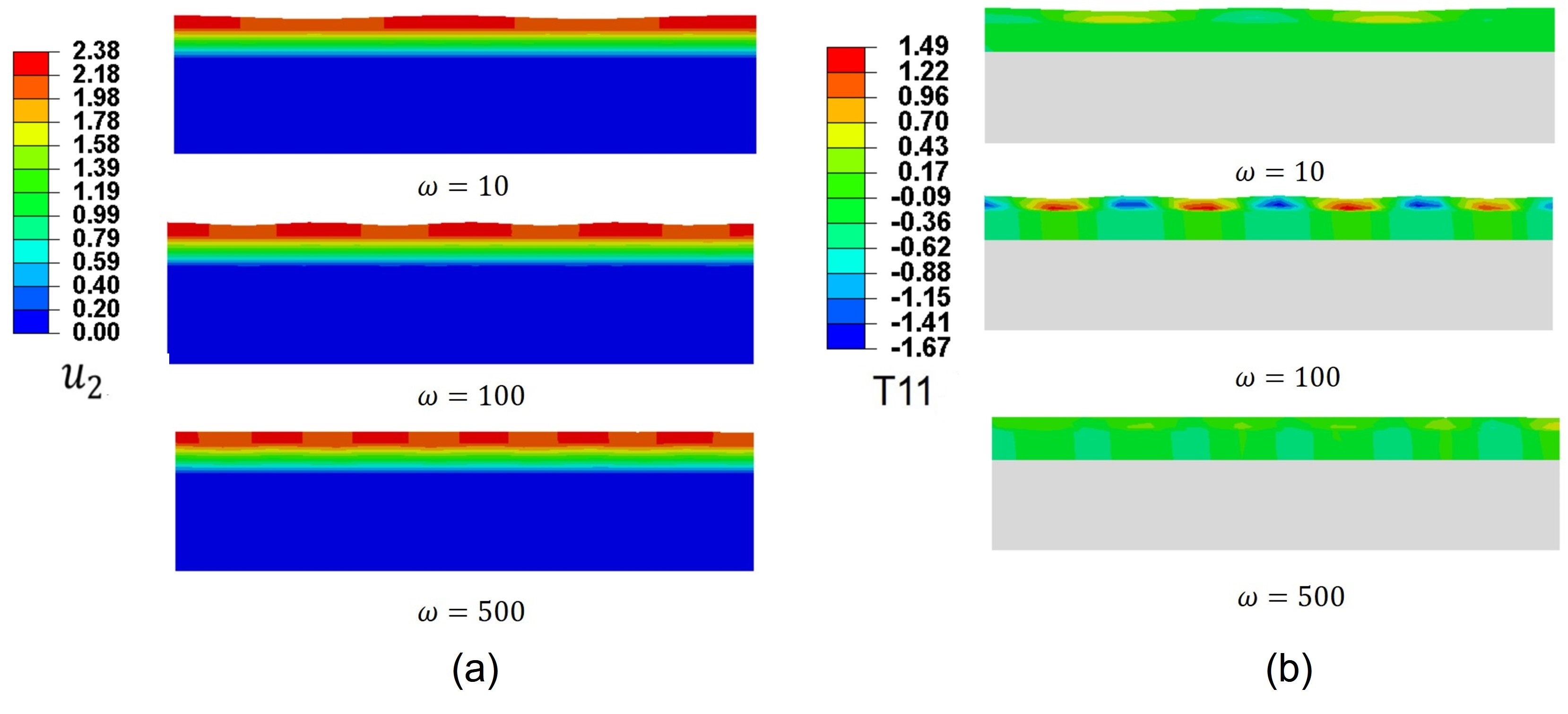} 
     \centering
     \end{tabular}   
  \caption{Contour plots of vertical displacement in $\mu m$(left) and bending stress, T11 in MPa (right) for perturbation wavelengths 10, 100, and 500~$/ \mu m$ when plated at $0.5~\textrm{mA}/{\textrm{cm}}^2$ for 1 hour.}
    \label{w100}
\end{figure}

Figure \ref{w100} presents contour plots for the three wavelengths, depicting (a) vertical displacement and (b) bending stress contours. For $\omega=10/\mu m$ and $\omega=100/\mu m$, instability propagates, and causes wrinkles on the top surface. However, reducing the surface imperfection further ($\omega=500/\mu m$), the instability does not propagate, as evidenced by the flat top surface. This is corroborated by the nearly uniform bending stress observed at the SEI layer, as shown in the bottom-most contour plot of Figure \ref{w100}(b). \par
%%%%%% result fig 3
%%%%%
\subsection{Investigations of the role of SEI mechanical properties}
\begin{figure}[h!]
  \centering
  \begin{tabular}{c}
     \includegraphics[width=0.95\textwidth]{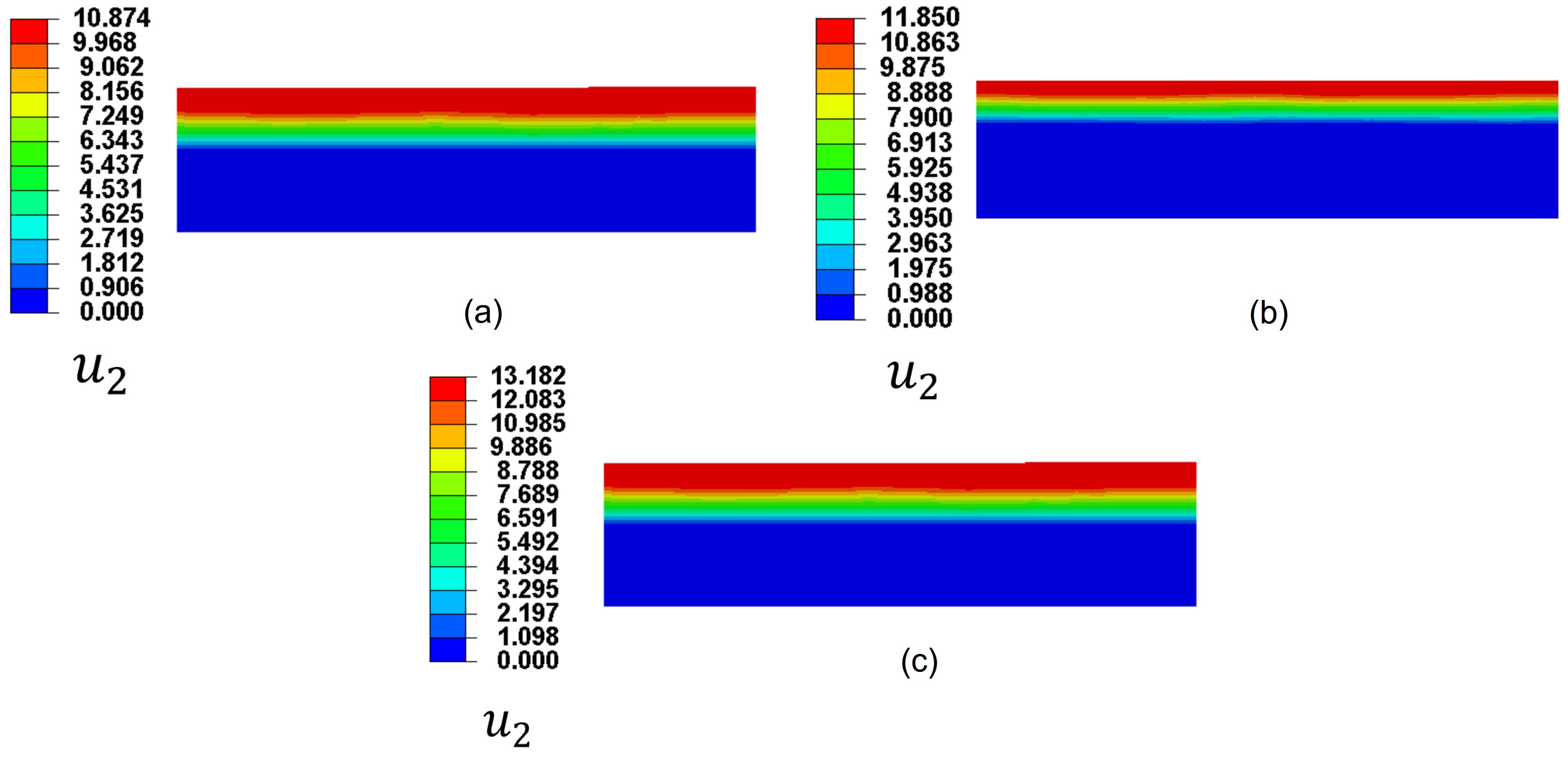}\\
     \centering
     \end{tabular}   
  \caption{Vertical displacement in $\mu m$ contour plots of: (a) SEI modulus 342 MPa, thickness $5~\mu m$, $\gamma=10~ N/m$; (b)  SEI modulus 30 GPa, thickness $5\mu m$; and, (c) SEI modulus 342 MPa, thickness $10~\mu m$ after 1 hour plated at $0.5~mA/{cm}^2$. }
    \label{Listab}
\end{figure}

The results obtained from the simulation presented in \fref{w100} demonstrate that the roughening of the anode depends on the wavelength of initial perturbation. To prevent the growth of roughening, it is necessary to keep it below certain values. However, in the real battery, it is impossible to exercise such control. To address the issue, mechanical approaches can be adopted, such as controlling the surface tension or the elastic modulus of the artificial SEI layer or keeping the SEI residual stress below a certain limit. To demonstrate this, we consider the following three cases: (i) a surface tension of $\gamma=10~N/m$ at the SEI/I interface with $E_{SEI}=342~MPa,\; \omega=10/\mu m, \; h_{\textrm{SEI}}=5~\mu m$; (ii) no surface tension with $E_{SEI}=30~GPa, \; h_{\textrm{SEI}}=5~\mu m$; and (iii) $E_{SEI}=342~MPa,\; h_{\textrm{SEI}}=10~\mu m$ with no surface tension. In all cases, $\omega=10/ \mu m$. In all cases, we observe a stable Li plating for 1 hour at a rate $0.5~ mA/{cm}^2$ as demonstrated by homogeneous displacement in both cases shown in \fref{Listab}.
These findings suggest that enhancing surface tension, elastic modulus, or SEI thickness, alone or in combination, enhances the stability of SEI and anode during Li-plating in a liquid electrolyte Li-metal battery. Using our new modeling and computational framework, we present some anode stability maps in the following section. For this purpose, we take the electrochemical properties of SEI as, $\kappa_{\textrm{SEI}}=1\times 10^{-6}~S/cm;~ \textrm D_{\textrm{SEI}}=1\times 10^{-8} \textrm{cm}^2/s$. 
%%%%%
\subsubsection{Interplay between SEI elastic modulus and surface tension}

\begin{figure}[h!]
  \centering
  \begin{tabular}{c}
     \includegraphics[width=0.65\textwidth]{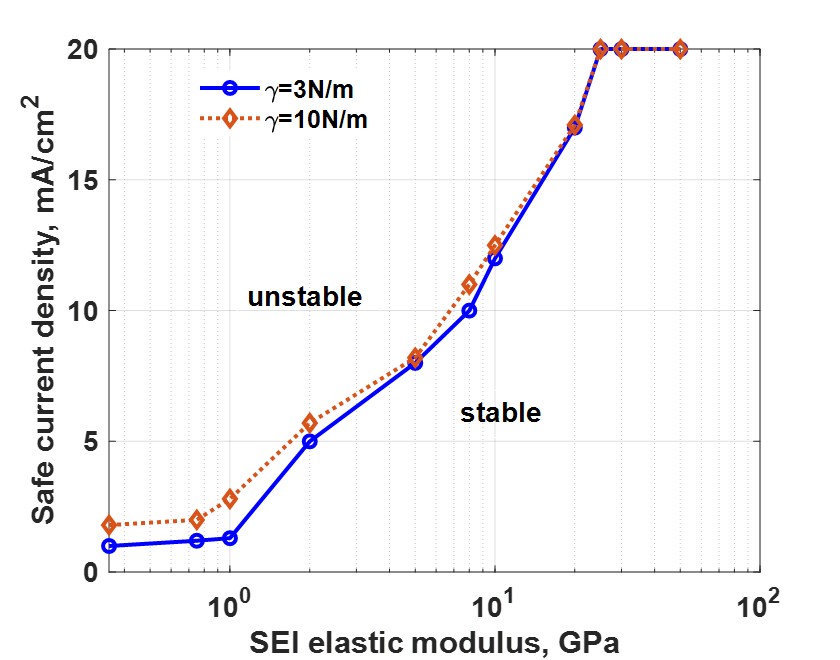}\\
     \centering
     \end{tabular}   
  \caption{Safe current density vs. SEI elastic modulus plot at 5$\mu m$ SEI thickness}
    \label{modulus}
\end{figure}

In this section, we explore the effects of SEI modulus and surface tension on anode stability. Numerical simulations are conducted with SEI elastic modulus ranging from 300~MPa to 100~GPa, assuming an initial perturbation wavelength of $\omega=10/\mu m$ and constant SEI thickness of $5~\mu m$. Plating is performed at various current densities until anode instability occurs within a maximum plating time of 10 hours. Our results indicate minimal influence of surface tension on stability below $3~ N/mm$, prompting an investigation into surface tension values ranging from $3-10 ~N/mm$. Figure~\ref{modulus} illustrates stable plating current density plotted against SEI modulus for surface tension values of $3~N/m$ (blue curve) and $10~N/m$ (red curve).
As observed from the plot, when the SEI is soft (modulus below 1~GPa), the safe current density is low, $\leq 2~mA/{cm}^2$. However, the 
\textit{safe current density} significantly improves at a SEI elastic modulus of 20~GPa or higher, $\geq 20~mA/{cm}^2$, and saturates above this point. Here, we define safe current density as the one where wrinkles do not appear upon Li-plating \par
As show in the plot, we observe that increasing the surface tension from $3~N/m$ to $10~N/m$ improves the SEI stability up to a SEI modulus of 10~GPa. Beyond this SEI modulus value, surface tension has negligible effect on increasing current density. As the modulus increasingly becomes higher, at a value greater than 20~GPa, surface tension plays no role in stability. This implies that softer SEI materials can increase plating current density by improving surface tension. However, for stiffer SEI, the instability suppression occurs via the generation of positive stress at the SEI layer rather than surface tension.
%%%%%%
%%%%%%
\subsubsection{Interplay between surface tension and SEI thickness} 

\begin{figure}[h!] 
  \centering
  \begin{tabular}{c}
     \includegraphics[width=0.65\textwidth]{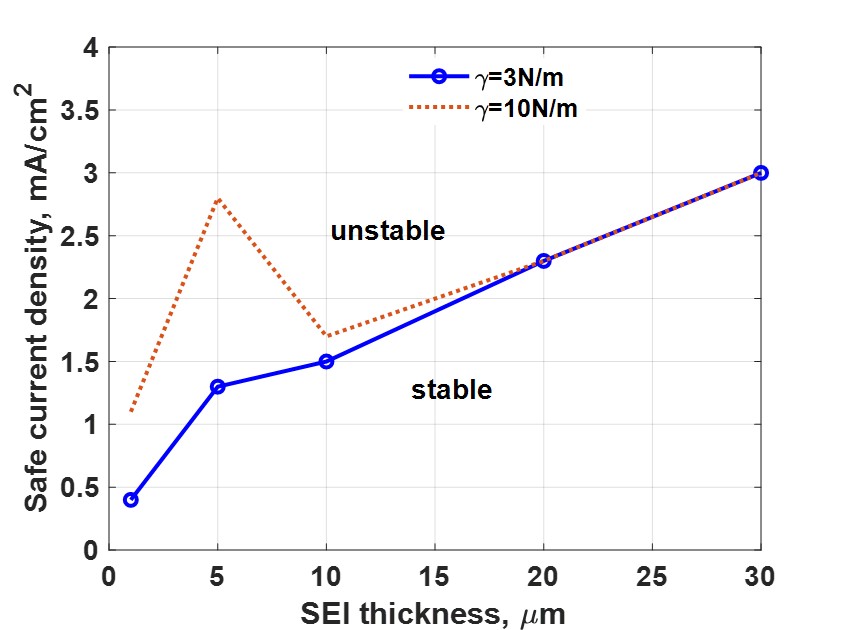}\\
     \centering
     \end{tabular}   
  \caption{ Safe current density vs. SEI thickness plot at 1GPa SEI elastic modulus}
  \label{thickness}
\end{figure}

In designing an artificial solid electrolyte interphase (SEI), it is important to ensure that the SEI has a finite thickness. This is because a positive mechanical blocking stress is generated, which is crucial in preventing dendrite formation. To investigate the impact of SEI thickness on anode stability, we conducted plating tests across different current densities, varying the thickness of the SEI from $1\mu m$ to $30\mu m$, assuming a moderate SEI elastic modulus of 1 GPa. Simulation results are shown in Figure~\ref{thickness}, with surface tension values of $3~N/m$ (blue) and $10~N/m$ (red). The plot shows that increasing the thickness of the SEI enhances stability by increasing the safe current density. However, it is worth noting that increased surface tension only boosts stability up to a thickness of $5~\mu m$. Beyond this thickness, the stable plating current density declines. This suggests that thicker SEI layers behave more like solid electrolytes, making surface tension enhancement ineffective in improving stability.
%%%%%
\subsubsection{Role of SEI residual stress} 
Compressive residual stress in the solid electrolyte interphase (SEI) layer is another major cause of mechanical instability in the anode of a battery. This stress can cause wrinkling of the anode and SEI layer, which creates hotspots for dendrite formation and growth \cite{liu2019wrinkling, wan2020mechanical, cho2020stress}. Although residual stress is a natural part of thin film deposition mechanics, it can be managed by engineering it to stay within certain limits. \par
\begin{figure}[h!] 
  \centering
  \begin{tabular}{c}
     \includegraphics[width=0.65\textwidth]{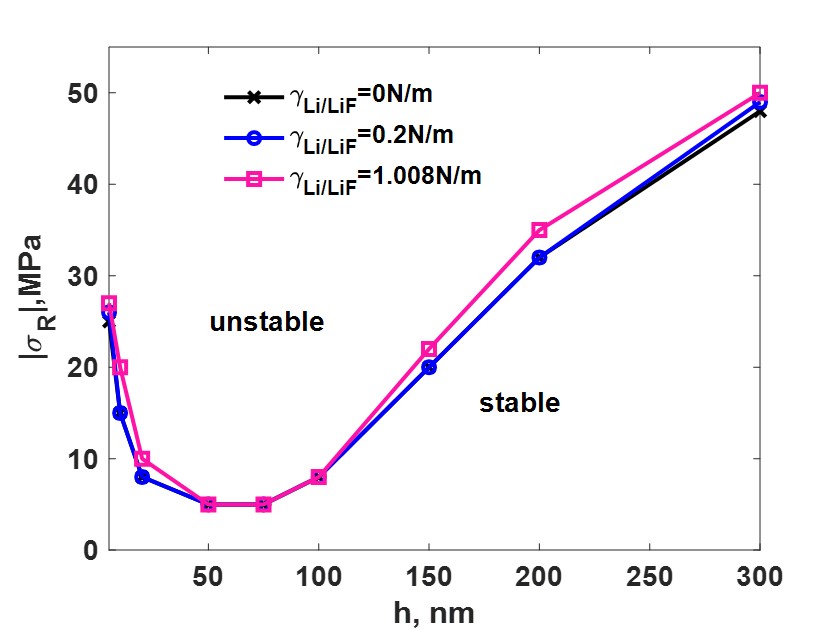}\\
     \centering
     \end{tabular}   
  \caption{Safe residual stress limit for LiF-type SEI for stable plating as a function of SEI thickness}
    \label{residual}
\end{figure}

We conducted an investigation to determine the safe residual stress limit for a LiF-type SEI with an elastic modulus of 65~GPa, for a thickness range of 5~nm to 300~nm at a low current density of $0.25 ~mA/{cm}^2$. The results are plotted over a period of 100 hours containing 5 plating and 5 stripping cycles. To further separate the role of residual stress only, here we neglect interfacial reaction, and consider the interphase grows proportionally to the initial incoming current. The results are plotted in \fref{residual}.\par
The results show that the stability of SEI (Solid Electrolyte Interface) generally increases as its thickness increases, except for an initial ultra-low thickness region. In this region, the thickness is almost zero, and therefore, the impact of residual stress is almost zero as well. However, beyond this point, increasing SEI thickness appears to be a better way to prevent SEI wrinkling. By controlling the residual stress at the SEI layer, we can avoid the unnecessary amplification of initial perturbation, thus minimizing anode roughening. Furthermore, we observed that surface tension has little effect on stability at SEI modulus 65~GPa. 

%%%%%%%
\subsection{Role of SEI inhomogeneities on stability}
\subsubsection{Modeling plating and stripping of Li with a SEI containing voids }

\begin{figure}[h!]
  \centering
  \begin{tabular}{c}
     \includegraphics[width=0.8\textwidth]{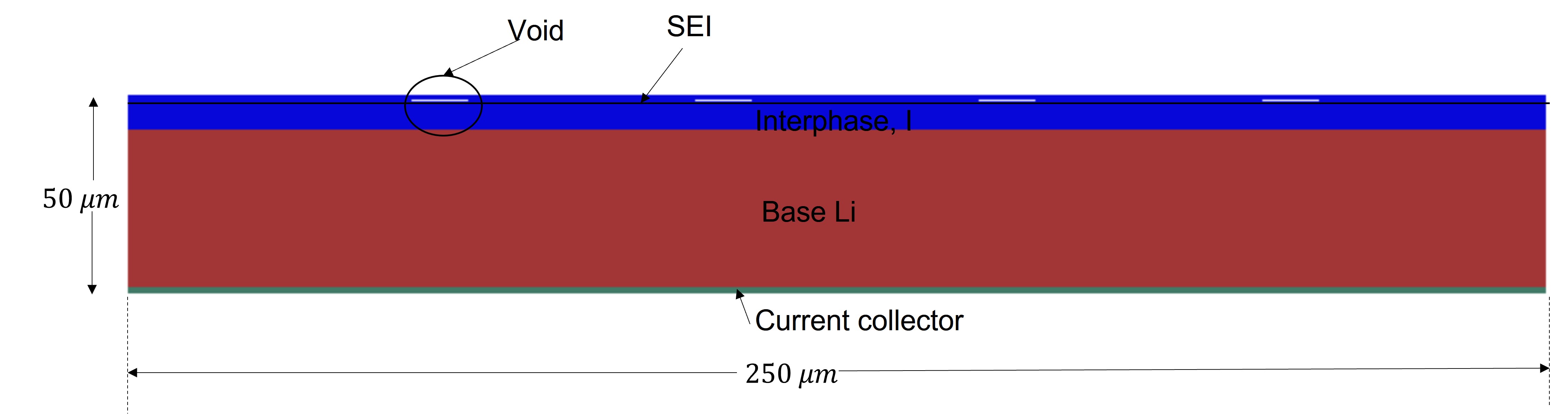}\\
     (a) \\
     \includegraphics[width=0.4\textwidth]{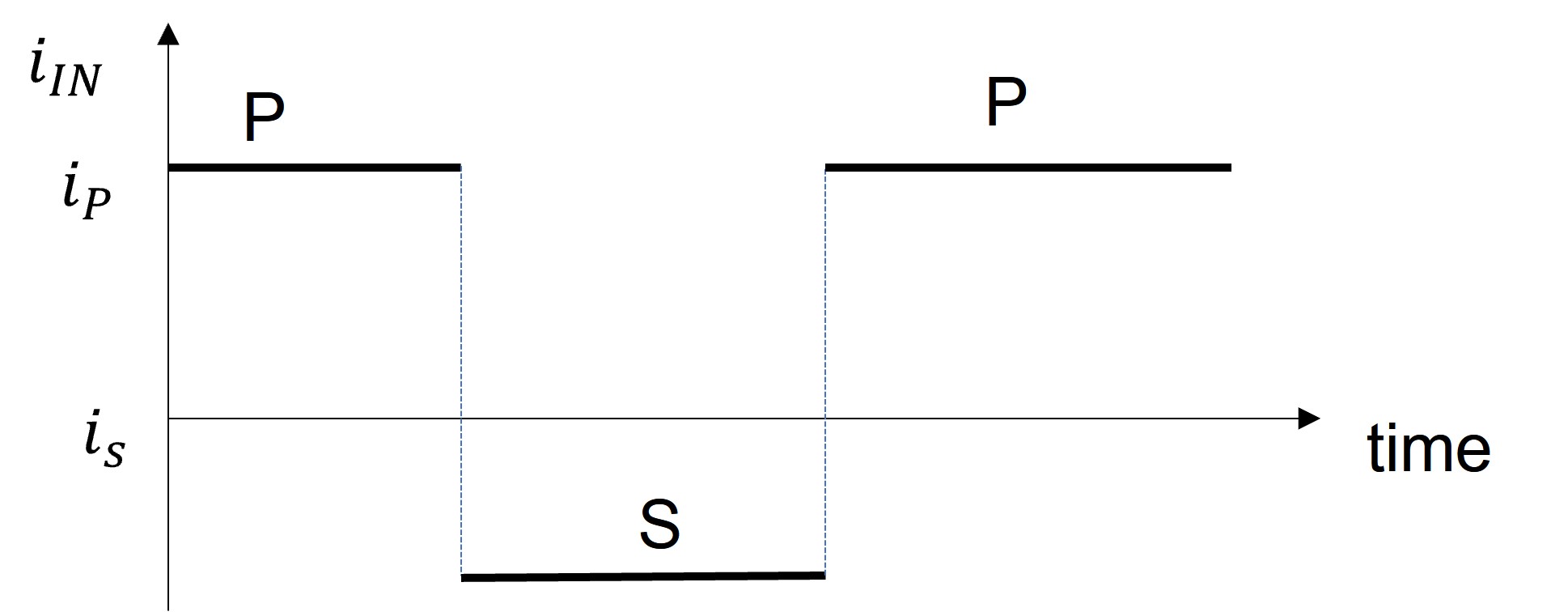}\\
     (b)\\
     \centering
     \end{tabular}   
  \caption{(a) Geometry of a Li half cell showing SEI voids; (b) galvanostatic plating and stripping boundary for the geometry in (a)  }
    \label{seivoid}
\end{figure}

The solid electrolyte interphase, both naturally formed or artificially sourced, often contains voids and impurities, which can cause voltage hotspots and stress inhomogeneities around these defects. Herein, we investigate the effect of the presence of voids at the SEI/I interface on anode instability during cyclic plating and stripping. For this purpose, we consider the geometry shown in \fref{seivoid}(a). Here, SEI layer has a thickness of 1~$\mu m$ containing 4 rectangular voids of dimension $1~\mu m \times 0.1 \mu m$ and located at the SEI/I interface. The interphase I is $5~\mu m$ thick. The properties of SEI are as follows: \par
E=340 MPa; $\nu$=0.342; $D_{\textrm{SEI}}=10~{\mu m}^2/s$, $\kappa_{\textrm{SEI}}=1e-4~ S/cm$, and $\Omega_{\textrm{SEI}}=2\times 10^{-4}~ \textrm{mol}/m^3$. The interphase and Li properties are taken from Table~\ref{table:MatParams} and \ref{table:Kinpar}. For the boundary condition, at first, we consider 1 hour of galvanostatic plating at a current density of $2~ mA/{cm}^2$, followed by 2020 seconds of striping at $2~ mA/{cm}^2$ and then another 3000 seconds of plating at $2~ mA/{cm}^2$, as shown in \fref{seivoid}(b). \par
%%%%%
\begin{figure}[h!]
  \centering
  \begin{tabular}{c}
     \includegraphics[width=0.97\textwidth]{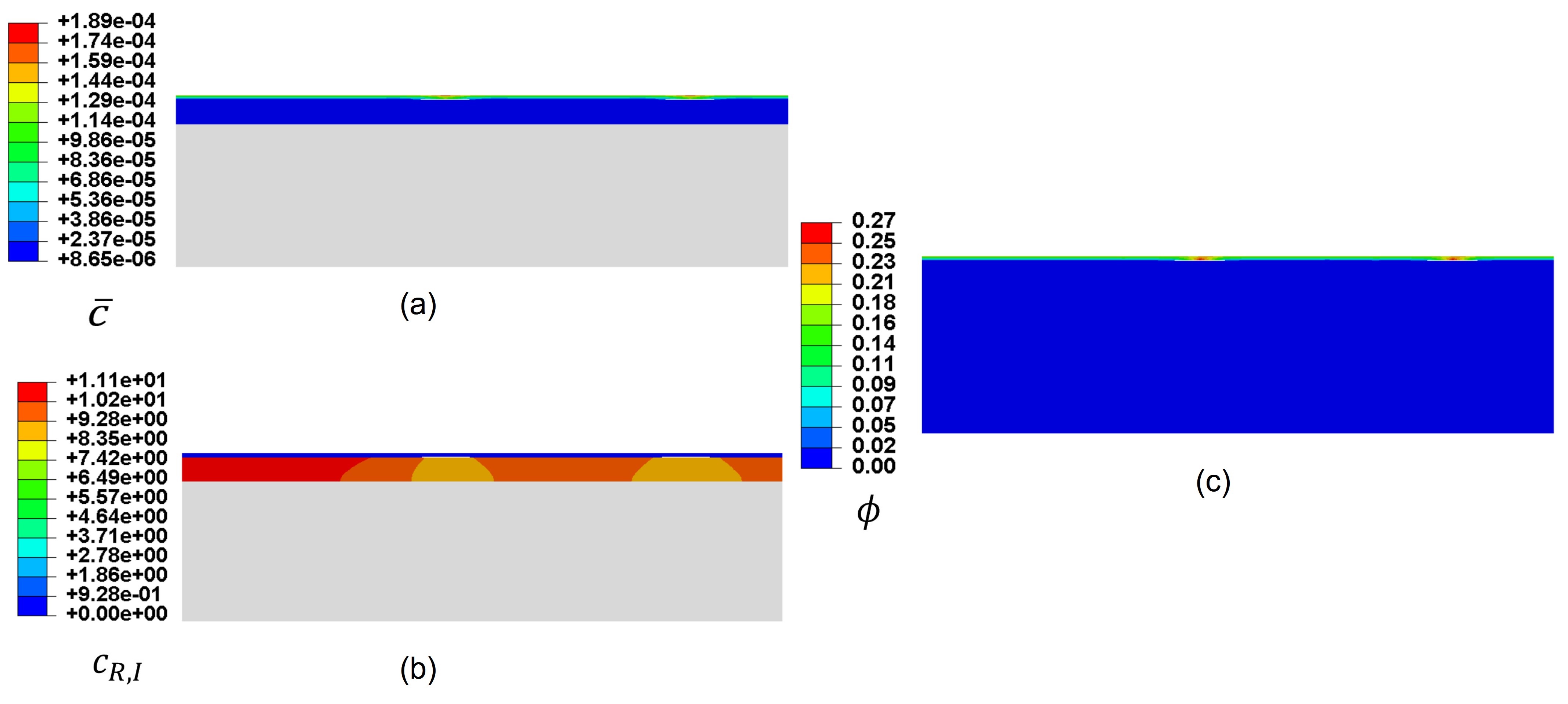}\\
     \centering
     \end{tabular}   
  \caption{Contour plots of (a) Normalized Li-ion concentration after 1s (b) Inhomogeneous Li-ion distribution at I after 1s (c) Voltage contour plot at the SEI and interphase I at the end of 1st plating cycle}
    \label{ion_voltage}
\end{figure}

%%%%%
The contour plots of normalized Li-ion concentration at the SEI layer and Li-ion concentration interphase I are shown in Figure~\ref{ion_voltage}(a) and (b), respectively. Due to the periodic nature of the voids, here we only show half of the contour plots. The presence of voids results in a higher ion concentration at the SEI around the voids, and the reaction current density around these voids is zero. Consequently, initially, Li-ion concentration is lower at the interphase I just below the voids. However, owing to the efficient transport properties of the SEI that we assume for this simulation, Li-ion homogenizes rapidly at the interphase (within about 4 seconds). The voltage contour is plotted in Figure~\ref{ion_voltage}(c). As expected, a much higher voltage is observed around the voids.\par
\begin{figure}[h!]
  \centering
  \begin{tabular}{c}
     \includegraphics[width=0.9\textwidth]{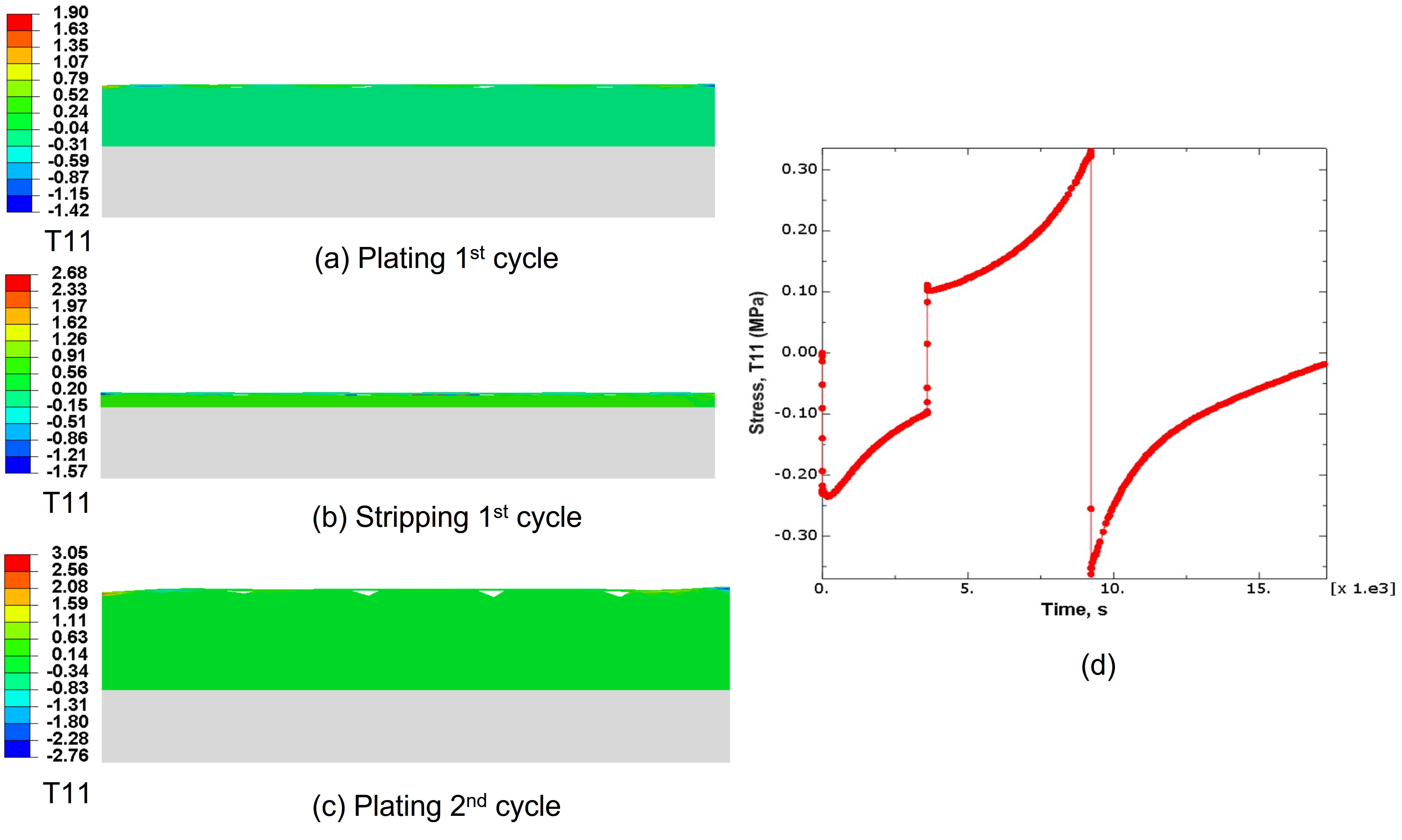}\\
     \includegraphics[width=0.9\textwidth]{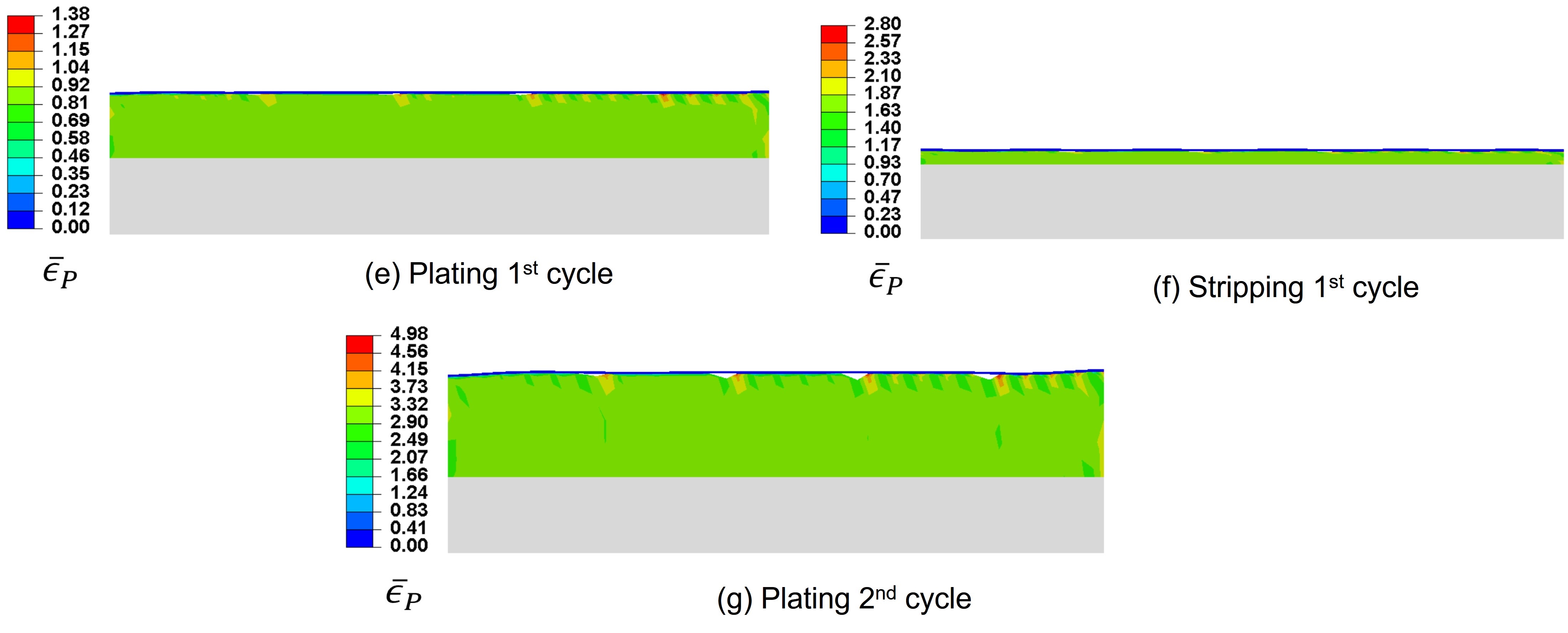}\\
     \centering
     \end{tabular}   
  \caption{(a)-(c) Contour plots of bending stress T11 for the SEI with voids geometry at the end of different cycles (d) Average bending stress in the SEI and interphase I; (e)-(f) Contour plots of equivalent plastic strain generated at the end of each cycle. }
    \label{T11_and_epsp_voids}
\end{figure}

%%%%%%
The contour plots of stress T11 shown at Figure~\ref{T11_and_epsp_voids}(a), (b) and (c) are at the end of 1st plating cycle, 1st stripping cycle and 2nd plating cycle, respectively. The average stress at the SEI and I layer is plotted in Figure~\ref{T11_and_epsp_voids}(d). Finally, Figure~\ref{T11_and_epsp_voids}(d)-(f) shows the contour plots of equivalent plastic strain at the end of each cycle. Analyzing the bending stress contours in Figure~\ref{T11_and_epsp_voids}(a-c), we observe that, stress consistently intensifies over the consecutive 3 plating and stripping cycles, indicating ongoing viscoplastic deformation in the interphase I. The voids deform and open after the 1st plating cycle, creating a slightly wavy surface. At this stage, the maximum plastic strain accumulates around the voids, as illustrated in Figure~\ref{T11_and_epsp_voids}(e). \par

Following the subsequent stripping cycle, the voids nearly disappear, and surface waviness reduces. However, a substantial amount of bending stress remains present in both the SEI and the interphase I layer, as shown in Figure~\ref{T11_and_epsp_voids}(b). A continued rise of plastic strain during stripping is observed in \fref{T11_and_epsp_voids}(f). In the subsequent plating cycle, the existing residual stress from the previous stripping cycle combines with the stress accumulated during plating, with high plastic strain (\fref{T11_and_epsp_voids}(c) and (g)); resulting in significant distortion of the voids. This, in turn, leads to pronounced roughening of the anode. Consequently, during continuous plating and stripping, the anode surface undergoes continuous evolution due to viscoplastic deformation of the interphase I, potentially causing instability even with improved electrochemical properties. \par
\fref{T11_and_epsp_voids}(d) shows that the average stress in the combined SEI+I assembly remains compressive during the second plating cycle. However, at the stripping cycle, stress becomes tensile due to the outflow of Li-ions from the interphase. As the final plating stage progresses, the combined stress gradually becomes tensile, which can lead to the fracture of the SEI upon further plating.\par
\begin{figure}[h!]
  \centering
  \begin{tabular}{c}
     \includegraphics[width=0.75\textwidth]{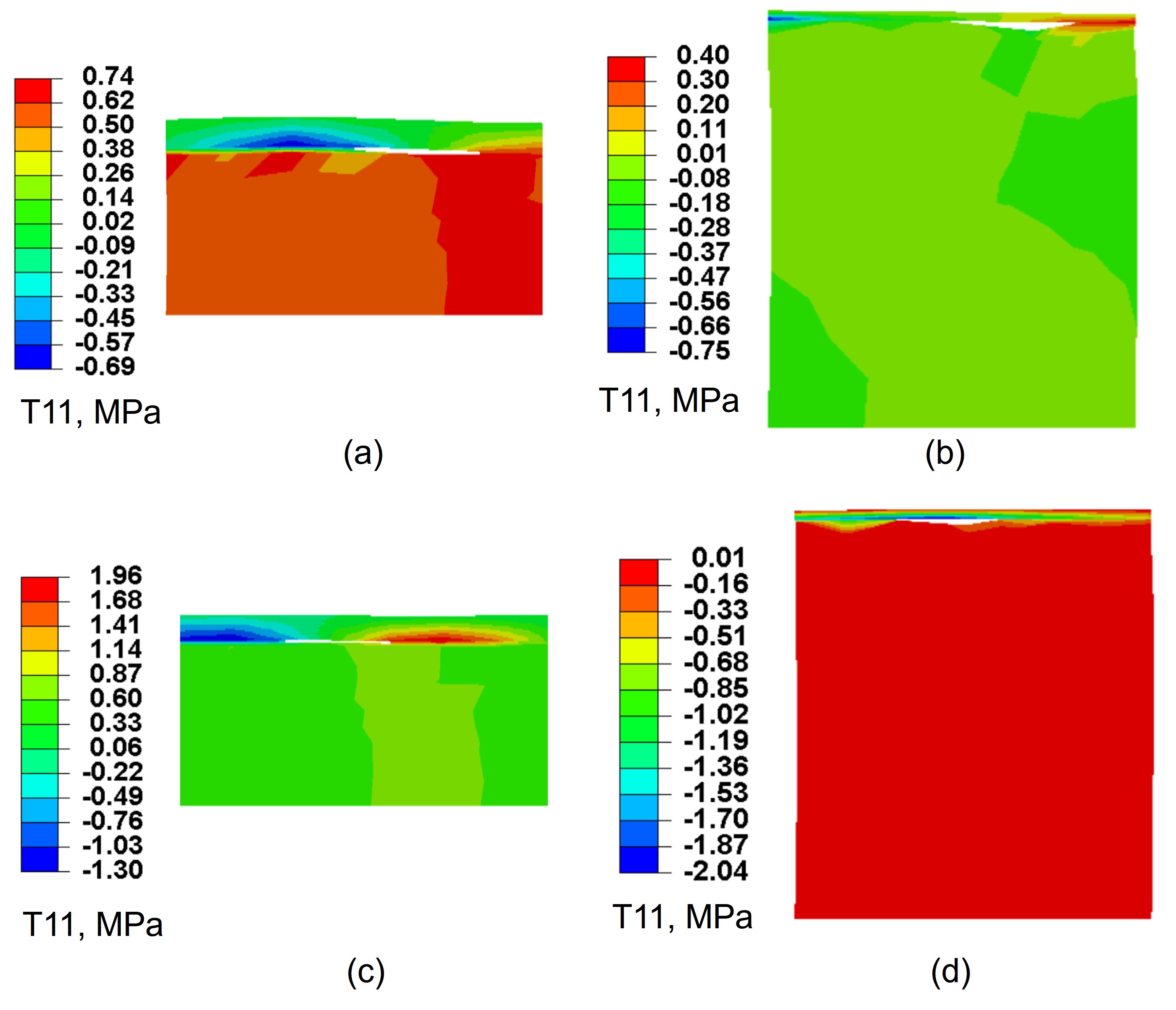}\\
     \centering
     \end{tabular}   
  \caption{Contour plots of stress distribution around a void with SEI elastic modulus of 340 MPa (a) after stripping cycle 1; (b) after plating cycle 2; and with SEI elastic modulus of 30 GPa (c) after stripping cycle 1; (d) after plating cycle 2.}
    \label{void_stress}
\end{figure}

%%%%%
Next, to explore the impact of SEI modulus in the presence of defects, we conducted a simulation with a SEI modulus of 30 GPa, in contrast to the previous run with a SEI modulus of 340 MPa. The stress distribution around a void is compared between these two SEI modulus settings, as illustrated in Figure \ref{void_stress}. A comparison of bending stress in the two modulus cases shows that the stiffer SEI exhibits significantly higher stress during the stripping cycle, attributed to its increased stiffness. Furthermore, during the last plating cycle, the stiffer SEI experiences a much higher magnitude of compressive stress. While the stiffer SEI appears capable of maintaining a flatter top surface compared to the softer SEI, the accumulated high stress around the void can lead to cracking due to its brittle nature. In contrast, the softer SEI seems more adept at accommodating larger local deformation for an extended period. This was evident in our simulation, where, during the 3rd cycle of plating, the simulation with the soft SEI ran for about 3000 seconds, while the stiff SEI case ran for 2000 seconds before abruptly losing convergence at the onset of instability. \par

The simulation results in this section indicate the impact of viscoplastic deformation and inhomogeneous stress distribution on the anode instability problem. The findings also emphasise the crucial need to address stress and electrochemistry as a coupled problem.

\subsubsection{Plating and stripping Li heterogeneous SEI proerties}
\begin{figure}[h!]
  \centering
  \begin{tabular}{c}
     \includegraphics[width=0.7\textwidth]{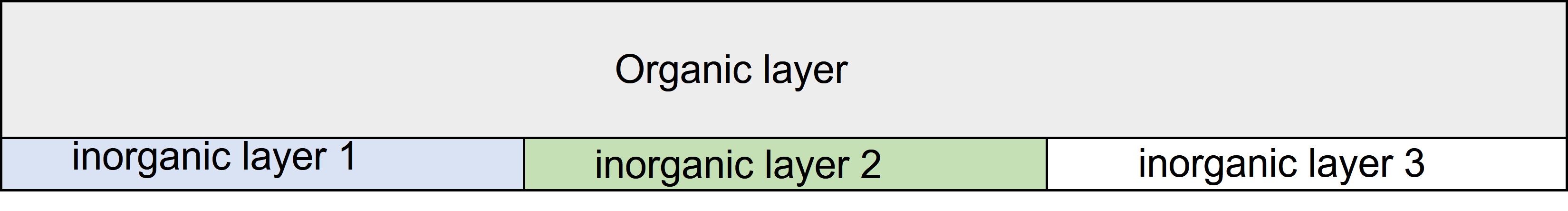}\\
     \centering
     \end{tabular}   
  \caption{Geometry of the composite SEI layer}
    \label{natualSEIgeo}
\end{figure}

To identify the role of SEI heterogeneous properties on stability, we conduct a numerical study of a composite SEI structure, mimicking the construction of natural SEI containing an organic layer followed by a heterogeneous inorganic layer composed of three inorganic layers, as shown in \fref{natualSEIgeo}. The upper layer of thickness $0.9 \mu m$, is a polymer layer, followed by a thin inorganic layer spanning $100 nm$. The properties of different SEI layers are outlined in Table~\ref{table:natsei}. \par 
%%%%
{\renewcommand{\arraystretch}{1}
\begin{table}[!h]
\caption{Properties of the multi-phase SEI layer}
\begin{center}
\begin{tabular}{|cccccc|}
\hline
\hline
Layer No & E(GPa) & $\nu$  & D ($m/s^2$)  & $\kappa$($S/cm$) & $\Omega_{\textrm{SEI}}$, $\textrm{mol}/m^3$\\
\hline
Organic layer & 0.35 & 0.3 & $1\times 10^{-12}$ & $1\times 10^{-6}$ & $2\times10^{-4}$\\
Inorganic layer 1 & 154 & 0.3 & $1\times 10^{-12}$ & $1\times 10^{-6}$ & $2\times10^{-4}$\\
Inorganic layer 2 & 108 & 0.3 & $1\times 10^{-12}$ & $1\times 10^{-6}$ & $2\times10^{-4}$\\
Inorganic layer 3 & 65 & 0.3 & $2\times 10^{-12}$ & $1\times 10^{-6}$ & $2\times10^{-4}$\\
\hline
\end{tabular}
\end{center}
\label{table:natsei}
\end{table}}
%%%%%%
We test 3 different cases for this study: (i)Two plating cycles at current density $0.5mA/{cm}^2$ for 1~hour with a 1~hour OCV hold in between; (ii) One plating cycle for 1 hour followed by a stripping cycle of 1 hour at current density $0.5 mA/{cm}^2$; (iii) plating at $2 mA/ {cm}^2$. \par
\begin{figure}[h!]
  \centering
  \begin{tabular}{c}
     \includegraphics[width=0.95\textwidth]{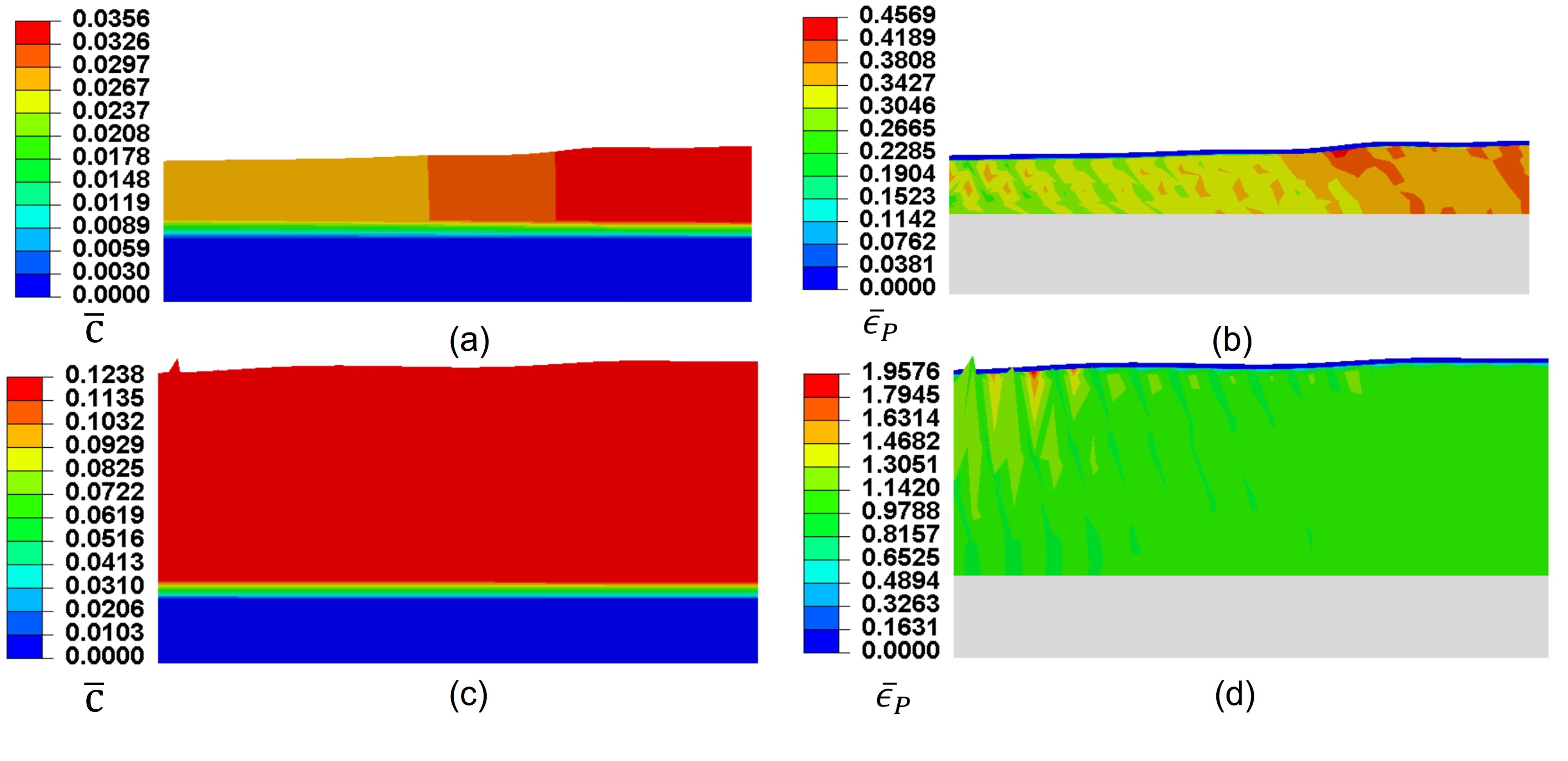}\\
     \centering
     \end{tabular}   
  \caption{Contour plots of (a) Normalized Li-ion concentration and (b) plastic strain after cycle-1; (c) Normalized Li-ion concentration and (d) plastic strain after cycle-2 for multi-layer composite SEI/I assembly when plated at current rate $0.5~mA/{cm}^2$. }
    \label{naturalplatingat0p5}
\end{figure}
In Figure \ref{naturalplatingat0p5}, we illustrate case (i) with a magnified y-scale to enhance visualization. The contour plots in Figure \ref{naturalplatingat0p5}(a) and (b) showcase the normalized concentration of Li-ion and plastic strain, respectively, after the initial cycle. Similarly, Figure \ref{naturalplatingat0p5}(c) and (d) reveal these plots after the second cycle. During the first cycle, the Li-ion concentration demonstrates significant non-uniformity, leading to greater interphase I growth beneath inorganic layer 3 compared to the other two layers. This phenomenon arises from the twofold higher diffusion coefficient of Li ions in inorganic layer-3, which dominates interphase growth through electrochemical properties. As a result, the plastic strain is also most pronounced in this region. At the end of the second plating cycle, the concentration of Li-ion becomes uniform throughout the entire solid electrolyte interface (SEI). However, the plastic strain reaches its maximum below the most rigid inorganic layer-1. This heightened plastic strain generates the highest equivalent tensile stress below layer-1, which causes severe distortion to the interphase elements. We can observe this phenomenon in (c) and (d) of \fref{naturalplatingat0p5}, where a spike-like protrusion is evident. This distortion could potentially cause the SEI to crack at this location, leading to dendrite growth. \par
\begin{figure}[h!]
  \centering
  \begin{tabular}{c}
     \includegraphics[width=0.98\textwidth]{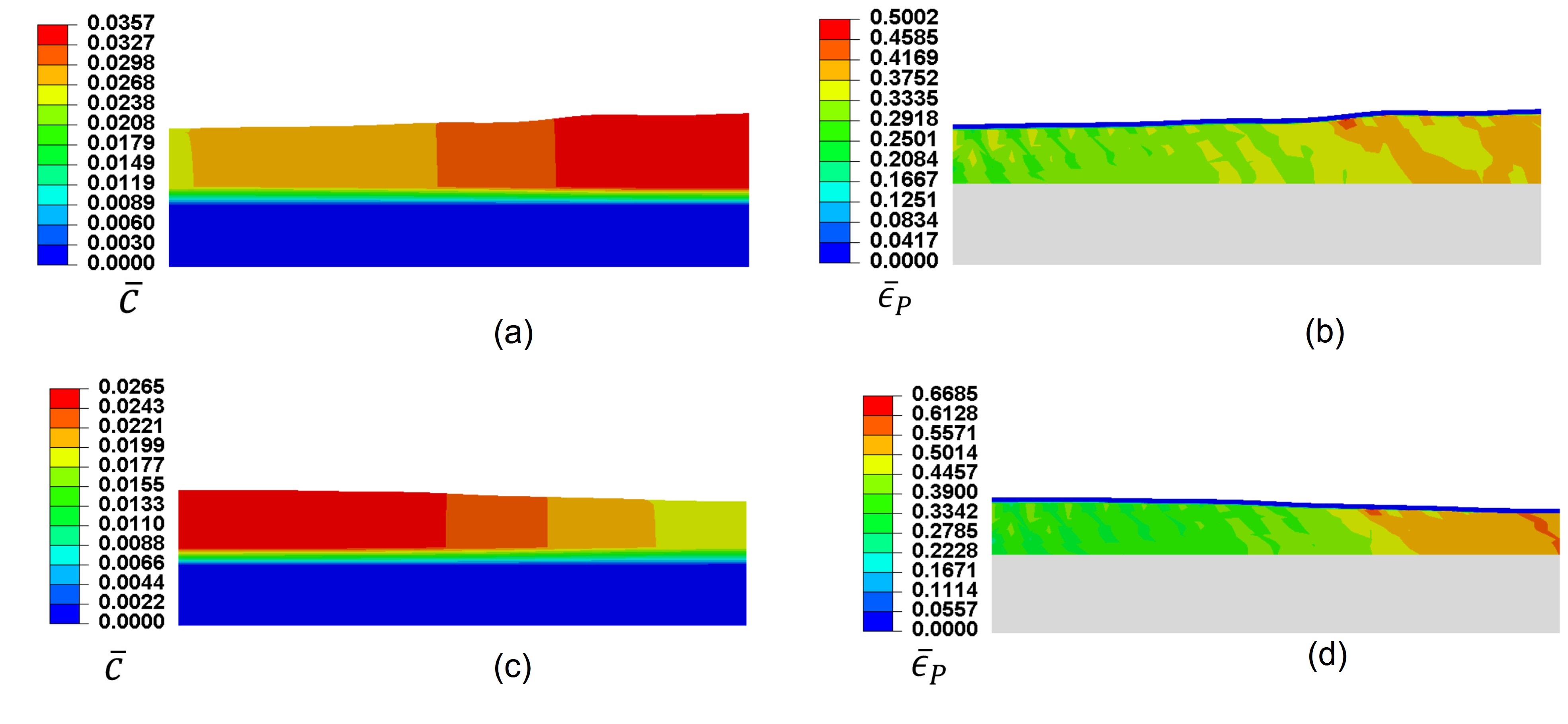}\\
     \centering
     \end{tabular}   
  \caption{Contour plots of (a) Normalized Li-ion concentration and (b) plastic strain after plating cycle; (c) Normalized Li-ion concentration and (d) plastic strain after stripping cycle for multi-layer composite SEI/I assembly when plated and stripped at current rate $0.5mA/{cm}^2$.}
    \label{naturalps}
\end{figure}

We have presented the results of case-(ii) in \fref{naturalps}, which involves a plating process followed by a stripping process. Similar to the previous scenario, interphase I experiences accelerated growth under the SEI inorganic layer-3 during plating. Due to the improved transport properties of this layer, we observe a similar trend during stripping, where Li-ion departure is faster in this region. As a result, the surface of the specimen after the stripping cycle is completely different from that after the plating cycle, with the highest height of protrusion found on the left side. It is worth noting that the plastic strain continues to accumulate primarily in the interphase elements beneath layer-3. This accumulation could potentially result in delamination of the SEI layer from the interphase layer due to ratcheting. \par
\begin{figure}[h!]
  \centering
  \begin{tabular}{c}
     \includegraphics[width=0.98\textwidth]{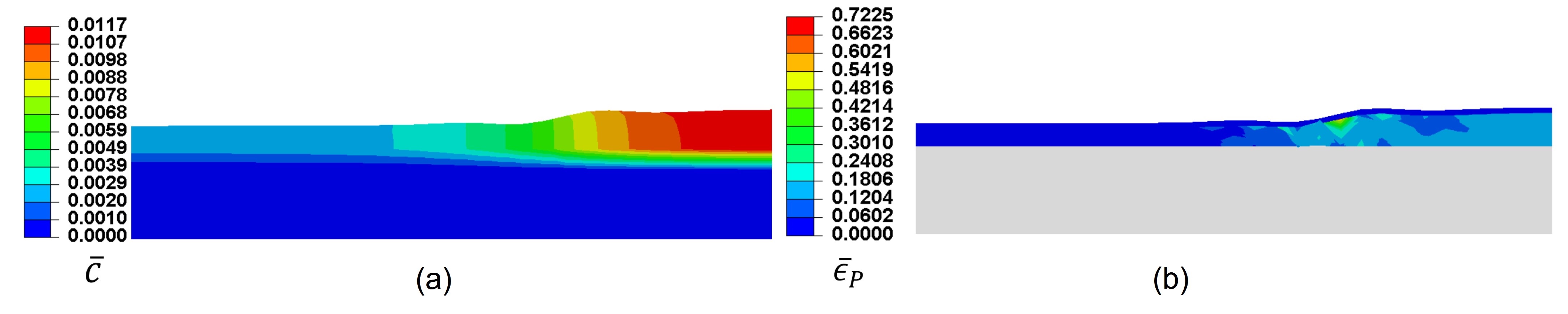}\\
     \centering
     \end{tabular}   
  \caption{Contour plots of (a) Normalized Li-ion concentration and (b) plastic strain for multi-layer composite SEI/I assembly when plated at current rate $2~mA/{cm}^2$.}
    \label{naturalh}
\end{figure}

Finally, in \fref{naturalh}, we observe the scenario where plating occurs at a high current density of $2~mA/{cm}^2$, corresponding to case-(iii). In this instance, protrusions rapidly emerge from SEI layer-3, occurring just after ~340s. The anode surface swiftly becomes unstable, reaching a point where it fails to accommodate the substantial deformation associated with further plating. 
\subsection{Double layer composite SEI stability}
\begin{figure}[h!]
  \centering
  \begin{tabular}{c}
     \includegraphics[width=0.9\textwidth]{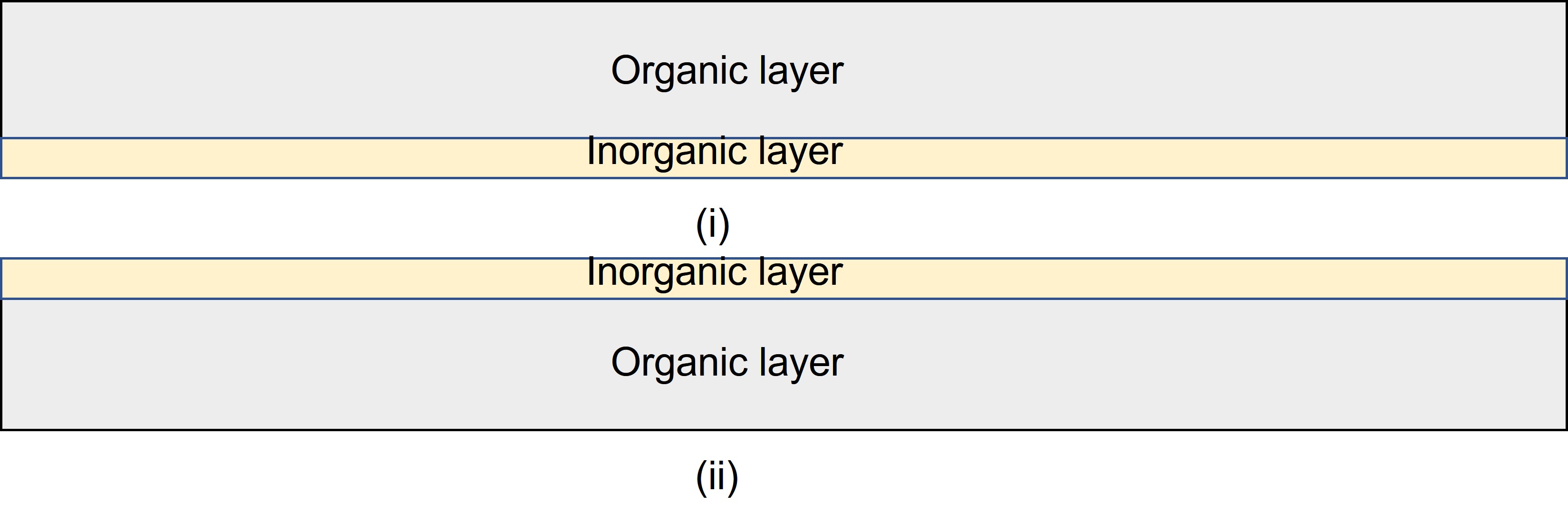}\\
     \centering
     \end{tabular}   
  \caption{Double layer composite SEI: (i) Organic layer on top of an inorganic layer; (ii) Inorganic layer on top.}
    \label{dlsei}
\end{figure}

The organic layer of the natural SEI enhances ductility and acts as a protective barrier between the anode and electrolyte. In conventional designs, artificial SEI is structured in a similar way, with the inorganic layer in direct contact with the anode and the organic layer placed above it. However, recent innovative designs have demonstrated that a configuration with an inorganic layer atop an organic layer performs more efficiently. This arrangement adheres better to the anode and provides greater flexibility to accommodate substantial volume changes during plating and stripping in Li-anode batteries \cite{liu2020review, liu2016artificial}. This section explores similar virtual experiments for both design approaches. The double-layer SEI, as depicted in \fref{dlsei}, is considered for this study. A polymer layer of $2.5~ \mu m$ thickness and a modulus of 340~MPa works as an organic layer while the inorganic layer is $0.5~ \mu m$ thick and has a modulus of 50~GPa. We plate Li at a current density of $1~ mA/{cm}^2$ for 1 hour for both types of SEI layers. A surface tension of $3~ N/m$ is assumed for both cases. For brevity, we refer to the configuration in \fref{dlsei}(i) as Polymer+Ceramic Double Layer (PCDL) and the design in \fref{dlsei}(ii) as Ceramic+Polymer Double Layer (CPDL).
\par

\begin{figure}[h!]
  \centering
  \begin{tabular}{c}
     \includegraphics[width=0.98\textwidth]{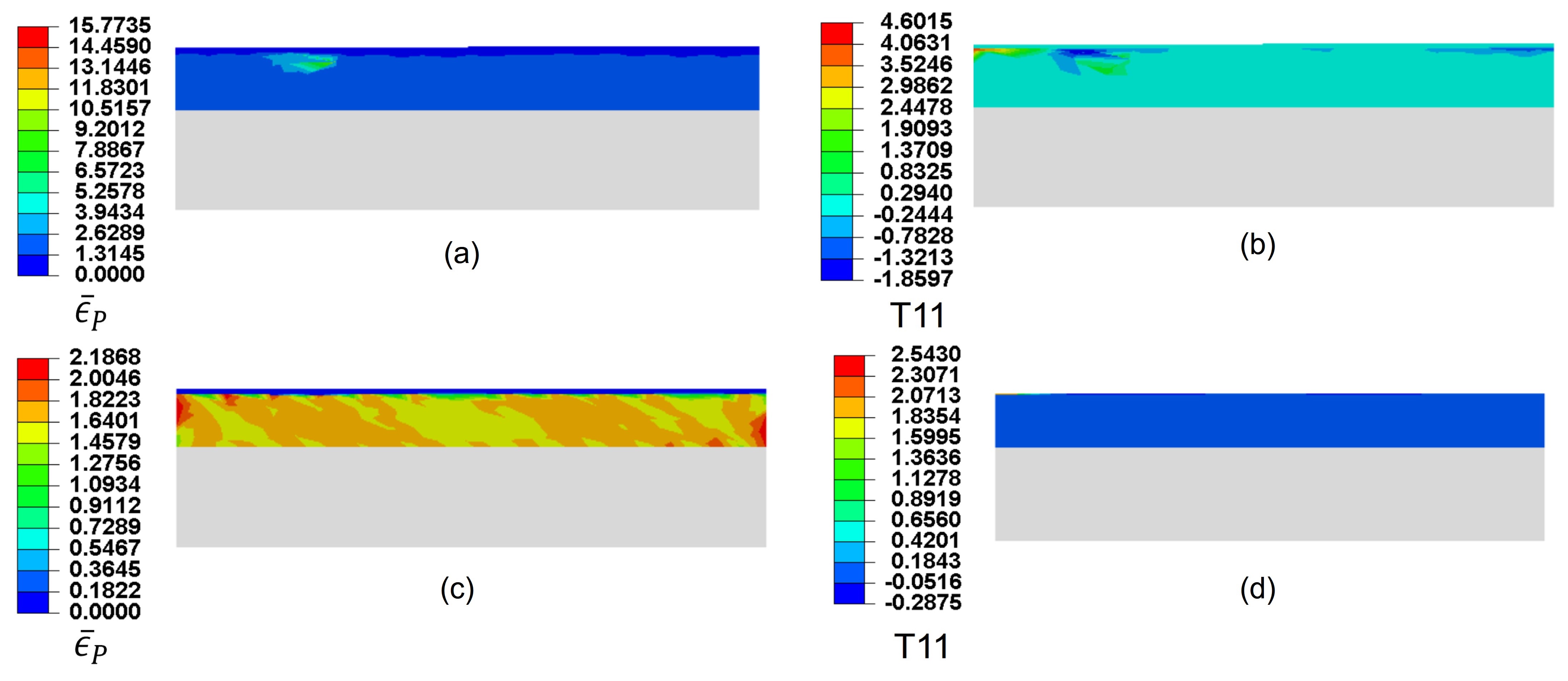}\\
     \centering
     \end{tabular}   
  \caption{Contour plots of (a) Equivalent plastic strain, and, (b) bending stress T11 for  polymer+ ceramic dl SEI; (c)Equivalent plastic strain, and (d) bending stress T11 for  ceramic+polymer dl SEI;}
    \label{dlcontour}
\end{figure}

In the scenario of the Polymer+Ceramic Double Layer (PCDL) SEI, the plating simulation concluded after a total runtime of 42.67 minutes, while, conversely, the Ceramic+Polymer Double Layer (CPDL) configuration sustained plating for up to 51 minutes. \fref{dlcontour} presents contour plots of bending stress and equivalent plastic strain for both cases after 42 minutes. For the PCDL SEI, localized high plastic strain is evident (\fref{dlcontour}(a)), whereas significantly lower plastic strain is observed in the CPDL layer. This difference arises from the stiffer layer in PCDL attempting to suppress non-local perturbations as the interphase grows, leading to localized plastic strain. In contrast, the CPDL configuration allows the interphase to deform freely without viscoplastic deformation, and the stiff top layer works to maintain a flat surface, resulting in a more homogeneous stress distribution, as depicted in \fref{dlcontour}(c). Additionally, stress T11 exhibits greater heterogeneity in PCDL (\fref{dlcontour}(b)), while it remains mostly homogeneous in the CPDL SEI (\fref{dlcontour}(d)). This suggests that a combination of hard-soft layers on the Li-anode provides better stability than the conventional soft-hard SEI layer by offering improved flexibility.
\section{Conclusion}
Our research presents new insights into anode instability. We developed a theoretical modeling and computational framework to study the effects of lithium's viscoplastic deformation with the interconnected dynamics of stress and electrochemistry in Li-anode. Contrary to Monroe and Newman's work, our findings suggest that high stiffness for SEI alone is not always stabilizing, rather a combination of SEI thickness and modulus is a better choice. We have demonstrated that the anode surface undergoes continuous evolution due to the viscoplastic behavior of Li, which makes the anode unstable after a finite number of charging and discharging cycles, due to its time and rate-dependent nature. We show that a bilayer SEI consisting of a soft inner layer in contact with Li-anode and a hard outer layer provides better stability. This is because the bilayer delays viscoplastic deformation of Li.
We studied the impact of artificial SEI properties and thickness on anode instability and obtained that an optimized artificial SEI with a carefully chosen thickness and modulus can enhance stability in liquid electrolyte Li-anode batteries by utilizing the positive effects of surface tension and elastic stress. Our study presents maps that illustrate the relationships between SEI modulus, surface tension, and thickness in shaping anode stability. These maps provide insights into the critical or safe plating current density. Despite these advancements, there is still much more to investigate and study. In particular, we need to better understand the mechanisms behind void formation during stripping, as well as the occurrences of anode and SEI cracking and delamination during repetitive plating and stripping. 

\section{Acknowledgements}
The authors thank the funding support from General Motors and Sayed Nagy, Caleb Reese, and Anil Sachdev from General Motors Global R{\&}D for valuable discussions.

\section*{Declaration of Interest}
The authors declare that they have no known competing financial interests or personal relationships that could have appeared to influence the work reported in this paper.

%\appendix

% References with bibTeX database %

\clearpage

%\section*{References}
%\bibliographystyle{elsarticle-num}
\bibliographystyle{elsarticle-harv}
%\bibliographystyle{elsarticle-num-names}
%\bibliography{biblist01,Master_Refs,ChesterRefs}
\bibliography{reference}

\end{document}